\documentclass[a4paper, 11pt]{article}

\usepackage[cp1251]{inputenc}
\usepackage[T2A]{fontenc}
\usepackage[english]{babel} 
\usepackage{amsfonts, amsmath, amssymb} 
\usepackage{graphicx}
\usepackage{enumerate}
\usepackage{flafter} 

\newtheorem{Th}{Theorem}
\newtheorem{Prop}{Proposition}
\newtheorem{Lm}{Lemma}
\newtheorem{Def}{Definition}
\newtheorem{Rm}{Remark}
\newtheorem{As}{Assumption}

\newcommand{\eps}{\varepsilon}
\newcommand{\R}{\mathbb R} 
\newcommand{\N}{\mathbb N} 
\newcommand{\Z}{\mathbb Z} 

\begin{document}

\title{On symplectic dynamics near \\a homoclinic orbit to 1-elliptic fixed point}
\author{L.\,Lerman, A.\,Markova\\
\normalsize Dept. of Diff. Equat. \& Math. Analysis and\\
\normalsize Research Inst. for Appl. Math. \& Cybernetics \\
\normalsize Lobachevsky State University of Nizhny Novgorod}

\maketitle

\begin{abstract}
We study the orbit behavior of a four dimensional smooth symplectic
diffeomorphism $f$
near a homoclinic orbit $\Gamma$ to an 1-elliptic fixed point under some natural genericity assumptions. 1-elliptic fixed point has
two real eigenvalues out unit circle and two on the unit circle. Thus
there is a smooth 2-dimensional center manifold $W^c$ where the restriction of
the diffeomorphism has the elliptic fixed point supposed to be generic (no strong resonances and first Birkhoff coefficient is
nonzero). Moser's theorem guarantees the existence of a positive measure
set of KAM invariant curves. $W^c$ itself is a normally hyperbolic
manifold in the whole phase space and due to Fenichel results
every point on $W^c$ has 1-dimensional stable and unstable smooth
invariant curves forming two smooth foliations. In particular, each KAM invariant curve has stable and unstable
smooth 2-dimensiona invariant manifolds being Lagrangian. The related stable and
unstable manifolds of $W^c$ are 3-dimensional smooth manifolds which are
assumed to be transverse along homoclinic orbit $\Gamma$. One of our
theorems present conditions under which each KAM invariant curve on $W^c$
in a sufficiently small neighborhood of $\Gamma$ has four transverse
homoclinic orbits. Another result ensures that under some Birkhoff
genericity assumption for the restriction of $f$ on $W^c$ saddle periodic
orbits in resonance zone also have homoclinic orbits though its
transversality or tangency cannot be verified directly.
\end{abstract}

\section{Introduction and set-up}

Any tools that can help to understand, if a given Hamiltonian system is integrable or non-integrable and therefore
has a complicated orbit behavior, are of the great importance. There are well known criteria based on the Melnikov
method \cite{Meln,Kozlov,LUm,Robinson,Delshams}, but they are mainly applicable for systems being nearly integrable.

There exists other class of criteria based on the study of the orbit behavior in definitely non-integrable systems:
if we know that some structures in the phase space are met only in non-integrable systems, then we may take the
existence of such a structure in the phase space of a system under consideration as a criterion of its non-integrability.
Such criteria are most efficient, if the structures mentioned can be rather easily identified. To this type of criteria
one can refer those based on the existence of homoclinic orbits to the different type of invariant sets, the most
popular are those related with homoclinic orbits to different types of equilibria, periodic orbits and invariant
tori \cite{Kozlov,Dev,Le1,Le2,KL1,KL2,KLDG,Bol,Cresson}.
Surely, non-integrability criteria are not the unique goal of the study: a much more interesting and hard problem
is to describe possible types of the orbit behavior in the system near such the structure and its changes when
parameters of the system vary.

In the paper we study a $C^r$-smooth, $r\ge 6,$ symplectic diffeomorphism $f$ on a $C^\infty$-smooth 4-dimensional
symplectic manifold $(M,\Omega)$, $\Omega$ is $C^\infty$-smooth non-degenerate 2-form. We assume $f$ to have an
1-elliptic fixed point $p$, that is, differential $Df_p$ has one pair of multipliers $e^{\pm i\alpha}$ on the unit
circle and a pair of real multipliers $\mu, \mu^{-1}$, $\mu\ne \pm 1$. Below we suppose $\mu$ to be positive and
$0< \mu < 1.$ We call such fixed point to be an orientable 1-elliptic point. The fixed point with  negative
$\mu, \mu^{-1}$ we call to be non-orientable. The non-orientable point becomes orientable, if one considers
$f^2$ instead of $f$.

Near an 1-elliptic fixed point there is a $C^{r-1}$-smooth 2-dimensional invariant symplectic center submanifold
$W^c$ corresponding to multipliers $e^{\pm i\alpha}$ \cite{Kelley,Mielke}. The restriction of $f$ on $W^c$ is a
$C^{r-1}$-smooth 2-dimensional symplectic diffeomorphism and $p$ is its elliptic fixed point. We assume $p$ to
be of the generic elliptic type \cite{AKN}, that is, strong resonances are absent in the system
$(\alpha \neq \displaystyle{\frac{\pi}{2},\frac{2\pi}{3}})$ and the first coefficient in the Birkhoff normal
form for $f \big|_{W^c}$ does not vanish. In this case we shall call an 1-elliptic fixed point to be
a \textit{generic 1-elliptic fixed point}. Then the Moser theorem \cite{Moser} is valid for the
restriction $f \big|_{W^c}$ near $p$, this gives a positive measure Cantor set of closed invariant curves on
$W^c$ which enclose $p$ and are accumulated to it. The needed minimal smoothness for a symplectic diffeomorphism
is 5 due to \cite{Russmann}. This explains the inequality $r\ge 6$.

Center manifold $W^c$ is a normally hyperbolic invariant manifold in the sense of \cite{Fenichel2,HPS} and has its local
$C^{r-1}$-smooth 3-dimensional stable manifold $W^{cs}_{loc}$ and local $C^{r-1}$-smooth 3-dimensional unstable
one $W^{cu}_{loc}$, since two other multipliers $\mu, \mu^{-1}$ are lesser than 1 and greater than 1, respectively
(these two local 3-dimensional manifolds for the fixed point $p$ are simultaneously center-stable and center-unstable manifolds,
respectively, this explains our notations). These manifolds can be extended till the global ones by the action of $f^{-1}$ and $f$, respectively.
The extended manifolds will be denoted as $W^{cs}$ and $W^{cu}$.

Each invariant KAM-curve $\gamma$ on $W^c$ can be considered as being saddle one, since it has local 2-dimensional
stable and unstable manifolds which can be also extended till global manifolds $W^s(\gamma)$, $W^u(\gamma)$ by the action of $f, f^{-1}$.
Topologically these manifolds are local cylinders, both being Lagrangian submanifolds in $M$ \cite{Arn}. The existence and smoothness of these
manifolds relies on the results of \cite{Fenichel1,Fenichel2} and will be proved in Appendix.

Fixed point $p$ has also two $C^r$-smooth local invariant curves through $p$ being its local stable $W^s_{loc}(p)$
and unstable $W^u_{loc}(p)$ manifolds \cite{Kelley}. Their extensions by the action of $f^{-1}$ and $f$ are $C^r$-smooth
invariant curves $W^s(p)$ and $W^u(p)$, respectively.

Our first two assumptions in the paper concern the existence of a homoclinic orbit to $p$ and its type.
\begin{As}[Homoclinic intersection]\label{intersection_as}
Curves $W^u(p)$ and $W^s(p)$ ha\-ve an intersection at some point $q$, thus generating a homoclinic orbit $\Gamma$ to fixed point $p$.
\end{As}

\begin{As}[Transversality condition]\label{transversality_as}
Manifolds $W^s(p)$ and $W^{cu}(p)$ are transverse at point $q$ and, hence, along $\Gamma$.
\end{As}

Later on in the section \ref{scattering_map_section} we will construct linear symplectic scattering map $S$ which
acts on tangent plane $T_pW^c$ and describes in the linear approximation an asymptotic behavior of orbits close
to $\Gamma$ after one-round travel near $\Gamma$. The restriction of differential $Df_p$ on symplectic invariant plane
$T_pW^c \subset T_p M$ is a linear symplectic 2-dimensional map with two eigenvalues $e^{\pm i\alpha}$,
and, therefore, this plane is foliated into closed invariant curves of the map. Every such a curve is an
ellipse, all of them can be obtained from the one multiplying their vectors at positive constants. Fix one
such ellipse $E \subset T_pW^c$. Then its image $S(E)$ is also an ellipse (usually not from the foliation) with the same center at the origin and
of the same area with respect to the restriction of 2-form $\Omega$ on this plane. Thus, the intersection
$E \cap S(E)$ consists of either four points (a generic case) or these two ellipses coincide (a degenerate case).
In the first case the intersection of two ellipses is transverse at every of four points.

\begin{As}[Genericity condition]\label{scattering_as}
The intersection $E \cap S(E)$ is transverse and therefore consists of four points.
\end{As}
It is evident that this assumption does not depend on the explicit choice of the ellipse $E$. This condition
allows one to select a generic case and provides the mean to verify this.

Our first result is the following theorem.
\begin{Th}\label{main_th}
Let a 4-dimensional symplectic diffeomorphism $f$ with 1-elliptic fixed point $p$ obeys Assumptions \ref{intersection_as},
\ref{transversality_as}, \ref{scattering_as}. Then there is a sufficiently small neighborhood $U$ of homoclinic orbit $\Gamma$
such that every closed invariant KAM-curve on $W^c(p)\cap U$ possesses four transverse homoclinic orbits in $U$.
\end{Th}

\begin{figure}
\centering
\parindent=0pt
\includegraphics[width=0.9\textwidth]{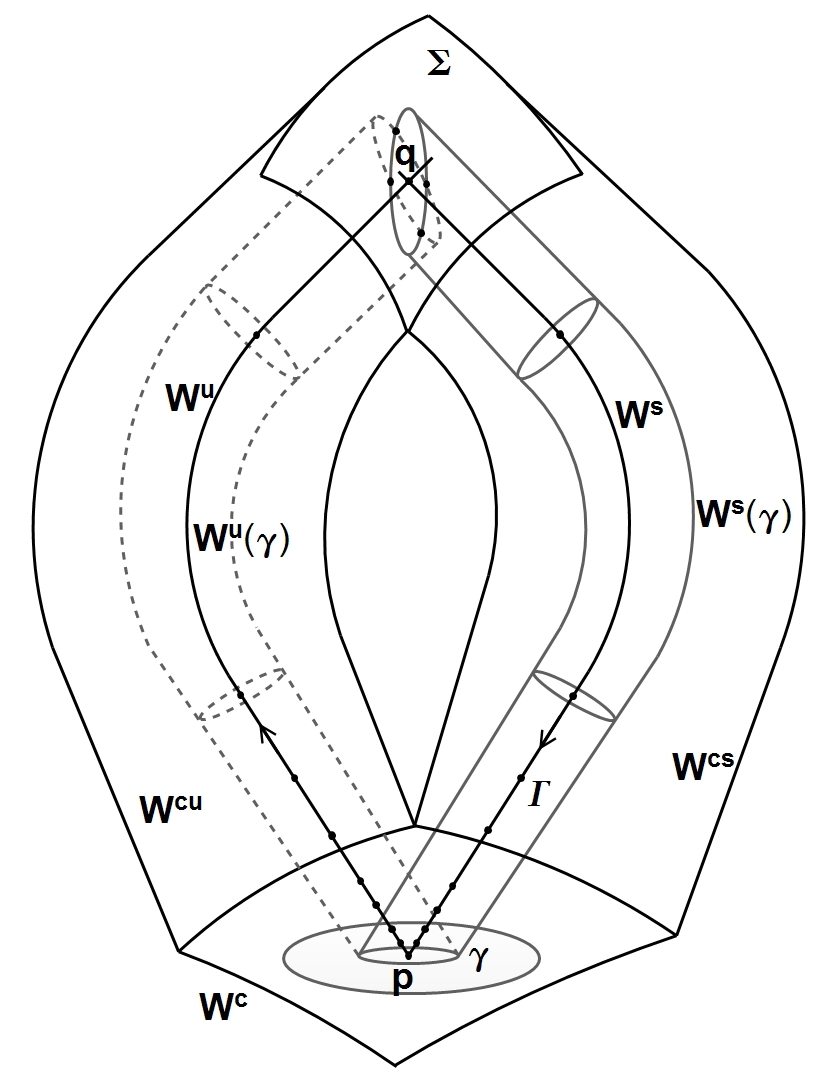}
\caption{Invariant manifolds in the neighborhood of homoclinic orbit}
\label{fig:invariant_manifolds}
\end{figure}

Intersection of invariant manifolds of the diffeomorphism $f$ in the neighborhood of homoclinic orbit are sketchy represented
on  Fig.~\ref{fig:invariant_manifolds}.
It is worth remarking that for our case center manifold $W^c$, as was
mentioned, is normally hyperbolic two-dimensional invariant manifold on which the
restriction of $f$ is a twist map. Thus our results on existence of
transverse homoclinic orbits to invariant KAM curves are connected with
the study of Hamiltonian dynamics near low-dimensional invariant whiskered
tori initiated in \cite{Easton} and extended in many recent papers (see, for instance, reviews \cite{dlL,DdlLS})

Before going to the proof, let us recall some related results for Hamiltonian vector fields
\cite{Le2,KL1,KL2,MS,GR1,GR2,MHO}. Homoclinic orbits to a saddle-center equilibrium for a real
analytic Hamiltonian system with two degrees of freedom, namely, for restricted circular three body problem,
were found numerically in \cite{LV} and proved to exist analytically through asymptotic expansions in \cite{MS}.
The problem on the orbit behavior of a real analytic Hamiltonian system near a homoclinic
orbit to a saddle-center equilibrium was first set up and partially solved in \cite{Le2}, though it was earlier discussed
in \cite{Conley}. In particular, the existence of four transverse homoclinic orbits to every small (Lyapunov's)
periodic orbit on the center manifold of the saddle-center was proved in \cite{Le2} using the Moser normal form
and the genericity condition was found first in \cite{Le2}. In \cite{GR1} under an additional assumption that a homoclinic orbit to
a saddle-center belongs to some invariant symplectic 2-dimensional submanifold (that is generically not the
case), the genericity condition was reformulated in terms of the related scattering problem for the transverse 2-dimensional system
linearized at the homoclinic orbit. It was first discovered in \cite{MHO} and in a more refined invariant form in \cite{GR2} that
in a generic 1-parameter unfolding of reversible 2 d.o.f. Hamiltonian systems that unfolds a Hamiltonian system
with a symmetric homoclinic orbit to a symmetric saddle-center equilibrium, there exists a (self-accumulated)
countable set of parameter values near the critical one such that for a point of this set the related Hamiltonian
system has a homoclinic orbit to its symmetric saddle-center. Usually these latter orbits are multi-round with respect to
the initial homoclinic orbit. Several applications, where non-integrability of a system under consideration
was proved using this method, can be found
in \cite{CN,GR3}. A partial extension of results to the case of Hamiltonian systems with $n$ degrees of freedom,
$n\ge 3$, having a center-saddle equilibrium (one pair of pure imaginary eigenvalues and the remaining ones with
nonzero real parts) with a homoclinic orbit, was given in \cite{KL2}. Here the scattering map was extended onto the
case when the center manifold is 2-dimensional but the dimension of transverse directions is $2(n-1)$.

In fact, the results we discuss here refer to a 3 d.o.f. Hamiltonian system on a $C^\infty$-smooth symplectic
manifold with a smooth Hamiltonian $H$ such that $X_H$ has a periodic orbit $\mathcal{C}$ of the center-saddle type. The
latter means the multipliers of this orbit (except for the common double unit) are a pair $e^{\pm i\alpha}$ and a pair of
reals $\mu, \mu^{-1},$ $\mu\ne \pm 1$. Such periodic orbit has 2-dimensional stable and unstable invariant
manifolds through $\mathcal{C}$, they both belong to 5-dimensional level $H=H(\mathcal{C})$. If these manifolds
have an intersection along some orbit $\Gamma$, then this homoclinic orbit tends to $\mathcal{C}$ as $t\to \pm \infty$.
Choose some cross-section $N$ to the flow through a point $p\in \mathcal{C}$ in 5-dimensional level $H=H(\mathcal{C})$.
We get a four-dimensional symplectic (w.r.t. the restriction of 2-form $\Omega$ to $N$) local Poincar\'e diffeomorphism
$f: N \to N$ with fixed point $p$ of the 1-elliptic type (corresponding to $\mathcal{C}$) defined in a neighborhood
$U \subset N$ of $p$. Intersection of stable and unstable manifolds of $\mathcal{C}$ with $N$ give smooth local curves
through the fixed point,
the traces of $\Gamma$ in $N$ form a countable set of homoclinic points accumulating at $p$. Fix one homoclinic point
$q_- \in U$ on the unstable curve and one homoclinic point $q_+ \in U$ on the stable curve. Choose some small
neighborhoods $V_- \subset U$ of $q_-$ and $V_+ \subset U$ of $q_+$ on $N$. Flow orbits define a symplectic map
$F:V_- \to V_+$, $F(q_-)=q_+$, that we call as global one. Then a symplectic first return map defined as $f$ for
points which belong to $U\setminus V_-$ and as $F$ for points in $V_-$ is a map we discuss.

The local center manifold $W^c(\mathcal{C})$ for periodic orbit $\mathcal{C}$ is of dimension four, it contains
the symplectic cylinder filled with periodic orbits (continuations of $\mathcal{C}$ onto close levels of $H$) and
if conditions of Theorem \ref{main_th} hold, then the restriction of the system on $W^c(\mathcal{C})$ has a positive
measure set of invariant 2-dimensional tori with Diophantine rotation numbers. When we fix the level $H=H(\mathcal{C})$,
then its intersection with the center manifold is 3-dimensional. Every torus $\mathcal T$ has stable and unstable
3-dimensional manifolds which intersect each other along four transverse homoclinic orbits to the torus within
5-dimensional level $H=H(\mathcal{C})$.

\section{Consequences of the transversality condition}\label{transversality_condition_consequence}

Due to  Assumptions \ref{intersection_as} and \ref{transversality_as}, two smooth 3-dimensional manifolds $W^{cs}(p)$
and $W^{cu}(p)$ intersect transversally at a homoclinic point $q$ and thus along a smooth 2-dimensional disk $\Sigma$
containing $q$. This disk is symplectic w.r.t. 2-form $\omega$ being the restriction of 2-form
$\Omega = dx \wedge dy + du \wedge dv$ on $\Sigma$. Indeed, in section \ref{homoclinic_orbits_existence_proof}
it will be proved that in normalized coordinates, in which $\Omega = dx\wedge dy + du\wedge dv$, disk (more exactly,
some its finite iteration under $f$) will have the following representation:
$$
x = \Phi(u,v),\; y = 0.
$$
This implies $\Sigma$ be symplectic w.r.t. 2-form $\omega = du \wedge dv$. The following lemma is valid:
\begin{Lm}\label{conv}
If Assumptions \ref{intersection_as} and \ref{transversality_as} hold, then $\Sigma$ and $W^{u}(p)$ are transverse
at $q$ within $W^{cu}(p)$ and, therefore, $W^u(p)$ and $W^{cs}(p)$ are also transverse at $q$.
\end{Lm}
{\bf Proof}. To prove this lemma, we use some symplectic coordinates $(x,u,y,v)$, $\Omega = dx \wedge dy + du \wedge dv$,
in a neighborhood $V$ of point $p$ in which manifolds $W^{cu}$ and $W^{cs}$ are straightened, that is they are given as
$x=0$ (for $W^{cu}$) and $y=0$ (for $W^{cs}$). In addition, in these coordinates local stable manifold $W^s(p)$ is
given as $y=u=v=0$ and local unstable manifold $W^u(p)$ is done as $x=u=v=0$. The existence of such coordinates is proved
in Appendix. We also assume that $q\in V$. Since orbit $\Gamma$ through $q$ is homoclinic, then there is an integer $N > 0$
such that $f^{-n}(q)\in V$ for all $n\ge N$. Denote $q_1$ the point $f^{-N}(q)\in V$ and let $l^u$ be the tangent to
$W^u(p)$ at $q_1$. Denote $L = Df_{q_1}^N: T_{q_1}M\to T_qM$, then $L(l^u)$ is transversal to $T_qW^{cs}$ in virtue to
Assumption \ref{transversality_as} (transversality condition). Set $D = L(T_{q_1}W^{cu})\cap T_qW^{cs}$, $D$ is 2-dimensional
plane. One needs to prove that $l^s$ (the tangent to $W^s(p)$ at $q$) does not belong to $D$, that is intersects $D$ at only one
point. For linear symplectic map $L$ the following matrix representation holds:
$$
L =
\begin{pmatrix}
	a & b \\ c & d \\
\end{pmatrix} =
\begin{pmatrix}
	a_{11} & a_{12} & b_{11} & b_{12} \\ a_{21} & a_{22} & b_{21} & b_{22} \\
	c_{11} & c_{12} & d_{11} & d_{12} \\ c_{21} & c_{22} & d_{21} & d_{22} \\
\end{pmatrix},
$$
where $a,b,c,d$ are $2\times 2$-matrices. Since $W^{cu}$, $W^{cs}$ are straightened  in coordinates we use, tangent
spaces to $W^{cu}$, $W^u(p)$ at $q_1$ and tangent spaces to $W^{cs}$, $W^s(p)$ at $q$ are written as follows:
$$
T_{q_1}W^{cu} =
\begin{pmatrix}
	0 \\ u \\ y \\ v \\
\end{pmatrix},\;
T_qW^{cs} =
\begin{pmatrix}
	\bar{x} \\ \bar{u} \\ 0 \\ \bar{v} \\
\end{pmatrix},\;
l^{u} =
\begin{pmatrix}
	0 \\ 0 \\ y \\ 0 \\
\end{pmatrix},\;
l^{s} =
\begin{pmatrix}
	\bar{x} \\ 0 \\ 0 \\ 0 \\
\end{pmatrix}.
$$
Transversality of $L(l^u)$ and $T_qW^{cs}$ is expressed as $d_{11}\ne 0$ in matrix $L$. Indeed, one has
$L(l^u) = (b_{11}y, b_{21}y, d_{11}y, d_{21}y)^T$ (vector-column). Transversality of $L(l^u)$ and $T_qW^{cs}$
means that determinant
$$
\begin{vmatrix}
	b_{11} & 1 & 0 & 0 \\
	b_{21} & 0 & 1 & 0 \\
	d_{11} & 0 & 0 & 0 \\
	d_{21} & 0 & 0 & 1 \\
\end{vmatrix} = d_{11}
$$
does not vanish.

The plane $D$ is given by the set of solutions of the system (\ref{tangent_plane}):
\begin{equation} \label{tangent_plane}
D:\left\{ {
\begin{matrix}
	\bar{x} = a_{12}u + b_{11}y + b_{12}v \\
	\bar{u} = a_{22}u + b_{21}y + b_{22}v \\
	0       = c_{12}u + d_{11}y + d_{12}v \\
	\bar{v} = c_{22}u + d_{21}y + d_{22}v \\
\end{matrix}
} \right.
\end{equation}
If $l^s \subset D$ then in the system above $\bar{u} = 0,\; \bar{v} = 0$ for all $(u, v)$. Expressing $y$
from the third equation in (\ref{tangent_plane}) and inserting into other equations we get a parametric
representation of plane $D$ (with parameters $u,v$). Consider separately subsystem (\ref{tangent_plane_subsystem}):
\begin{equation} \label{tangent_plane_subsystem}
\left\{ {
\begin{matrix}
	a_{22}u + b_{21}y + b_{22}v = \bar{u} \\
	c_{22}u + d_{21}y + d_{22}v = \bar{v}. \\
\end{matrix}
} \right.
\end{equation}
Due to inequality $d_{11} \neq 0$ we can express $y = (-c_{12}u - d_{12}v) / d_{11}$ from (\ref{tangent_plane})
and insert it into (\ref{tangent_plane_subsystem}):
\begin{equation} \label{tangent_plane_subsystem_1}
\left\{ {
\begin{matrix}
	d_{11}a_{22}u + b_{21}(-c_{12}u - d_{12}v) + d_{11}b_{22}v = \bar{u} \\
	d_{11}c_{22}u + d_{21}(-c_{12}u - d_{12}v) + d_{11}d_{22}v = \bar{v}. \\
\end{matrix}
} \right.
\end{equation}
Let us calculate the determinant of the system (\ref{tangent_plane_subsystem_1}). To this end, we rewrite
it in the following form:
$$
\left\{ {
\begin{matrix}
	u(d_{11}a_{22} - b_{21}c_{12}) + v(b_{22}d_{11} - b_{21}d_{12}) = \bar{u} \\
	u(d_{11}c_{22} - d_{21}c_{12}) + v(d_{22}d_{11} - d_{21}d_{12}) = \bar{v}. \\
\end{matrix}
} \right.
$$
This determinant is calculated as follows:
\begin{multline*}
\Delta = (d_{11}a_{22} - b_{21}c_{12})(d_{22}d_{11} - d_{21}d_{12})
            - (b_{22}d_{11} - b_{21}d_{12})(d_{11}c_{22} - d_{21}c_{12}) \\
= d^2_{11}a_{22}d_{22} - d_{11}a_{22}d_{21}d_{12} - b_{21}c_{12}d_{22}d_{11}
            + b_{21}c_{12}d_{21}d_{12} - d^2_{11}c_{22}b_{22} \\
+ b_{22}d_{11}d_{21}c_{12} + b_{21}d_{12}d_{11}c_{22} - b_{21}d_{12}c_{12}d_{21} \\
=d^2_{11}a_{22}d_{22} - d_{11}a_{22}d_{21}d_{12} - b_{21}c_{12}d_{22}d_{11} - d^2_{11}c_{22}b_{22} \\
+ b_{22}d_{11}d_{21}c_{12} + b_{21}d_{12}d_{11}c_{22} \\
= d^2_{11}(a_{22}d_{22} - c_{22}b_{22}) + d_{11}c_{12}(b_{22}d_{21} - b_{21}d_{22})
            + d_{12}d_{11}(b_{21}c_{22} - a_{22}d_{21}).
\end{multline*}
Matrix $L$ is symplectic, therefore the following identities hold (see, for instance \cite{GTSh}):
$$
a^T c \equiv c^T a,\; b^T d \equiv d^T b, \;d^T a - b^T c \equiv E =
\begin{pmatrix}
1 & 0 \\ 0 & 1 \\
\end{pmatrix}.
$$
The first identity is equivalent to equality:
$$
a_{12}c_{11} + a_{22}c_{21} = c_{12}a_{11} + c_{22}a_{21}.
$$
Similarly, the second matrix identity is reduced to equality:
\begin{equation} \label{symplectic_identity_2}
b_{12}d_{11} + b_{22}d_{21} = d_{12}b_{11} + d_{22}b_{21}.
\end{equation}
The third matrix identity gives the following relations:
\begin{equation} \label{symplectic_identity_3}
\begin{array}{ccc}
d_{11}a_{11} + d_{21}a_{21} - b_{11}c_{11} - b_{21}c_{21} & = & 1 \\
d_{11}a_{12} + d_{21}a_{22} - b_{11}c_{12} - b_{21}c_{22} & = & 0 \\
d_{12}a_{11} + d_{22}a_{21} - b_{12}c_{11} - b_{22}c_{21} & = & 0 \\
d_{12}a_{12} + d_{22}a_{22} - b_{12}c_{12} - b_{22}c_{22} & = & 1 \\
\end{array}
\end{equation}
Now, taking into account relations (\ref{symplectic_identity_2}), the second and the fourth equalities
in (\ref{symplectic_identity_3}), the expression for $\Delta$ can be transformed as follows:
\begin{multline*}
\Delta = d^2_{11}(1 + b_{12}c_{12} - d_{12}a_{12}) + d_{11}c_{12}(b_{11}d_{12} - b_{12}d_{11})
         + d_{12}d_{11}(d_{11}a_{12} - b_{11}c_{12}) \\
= d^2_{11} + d^2_{11}b_{12}c_{12} - d^2_{11}d_{12}a_{12} + d_{11}c_{12}b_{11}d_{12}
         - d^2_{11}c_{12}b_{12} + d^2_{11}d_{12}a_{12}\\
- d_{12}d_{11}b_{11}c_{12} = d^2_{11} \neq 0.
\end{multline*}
Thus, linear system (\ref{tangent_plane_subsystem_1}) has a unique solution $(u,v)$ at the given $(\bar{u},\bar{v})$.
So, $\bar{u}=0,\; \bar{v}=0$ only if $(u,v)=(0,0)$ and $l^s$ intersects $D$ at the unique point. $\blacksquare$

The Assumption \ref{intersection_as} says that $f$ is degenerate since generically two smooth curves in
a 4-dimensional manifold do not intersect. This assumption selects a codimension 2 set of diffeomorphisms
in the space of all $C^r$-smooth symplectic diffeomorphisms on $M$. Indeed, when a diffeomorphism with
a homoclinic orbit to an 1-elliptic fixed point is perturbed within the class of smooth symplectic ones,
for a perturbed $f'$ the fixed point $p'$ persists and its type is preserved. Therefore, due to transversality
condition, the intersection of perturbed $W^{cu}(p')$ and $W^{s}(p')$ persists as well, but the intersection
point does not give generically a homoclinic orbit to $p':$ with backward iterations of $f'$ the orbit
through the intersection point can be either a heteroclinic  orbit connecting $p'$ and some invariant curve
on $W^c(p')$ or some other orbit wandering near $W^c$ (recall that there are instability regions on $W^c(p')$,
the orbit returns to $W^c(p')$ staying within 3-dimensional $W^{cu}(p')$, thus it is locked between unstable
2-dimensional manifolds of invariant curves on $W^c(p')$, since they locally divide $W^{cu}(p')$).

Nevertheless, if we turn to the related 3 d.o.f. Hamiltonian system with a periodic orbit $\mathcal C$ of
1-elliptic type (or it can be called to be of the saddle-center type), then such an orbit belongs to a
smooth symplectic cylinder of periodic orbits of the same type. So, if $\mathcal C$ has a homoclinic
orbit, then for the related close levels of Hamiltonian on the cross-section to $\mathcal C$ one gets a
one-parameter family of symplectic Poincar\'e maps. Thus, if Hamiltonian itself depends on a parameter
$\eps$ in a generic way, then first return map for $\mathcal C$, derived by a homoclinic orbits to it, unfolds
to a two-parameter family of symplectic maps and hence any close smooth 1-parameter family of smooth Hamiltonians
also has a 1-elliptic periodic orbit with a homoclinic orbits to it. Thus, this phenomenon is generic for
generic 1-parameter unifoldings of a Hamiltonian with such the structure.

Now we return to the problem under study. In a neighborhood of homoclinic point $q$ let us consider
2-dimensional symplectic disk $\Sigma$ through $q$ being  the transverse intersection of extended
3-dimensional center-unstable manifold $W^{cu}(p)$ with 3-dimensional center-stable manifold $W^{cs}(p)$.
Below we shall prove the existence of smooth stable and unstable manifolds for any KAM-curve on $W^c$
lying in a sufficiently small neighborhood of $p$. All stable manifolds belong to $W^{cs}$ and all
unstable manifolds belong to $W^{cu}$. Hence, they intersect with $\Sigma$. The first statement
concerning this intersection is the following:

\begin{Lm}\label{equal}
Disk $\Sigma$ contains two Cantor sets of smooth closed curves $w_u(\gamma)$ and
$w_s(\gamma)$ being, respectively, traces of the related stable $W^s(\gamma)$ and unstable $W^u(\gamma)$
manifolds of invariant KAM-curves $\gamma \in W^c$. For a fixed invariant curve $\gamma$ integrals
of 2-form $\omega$ over disks $D^u(\gamma)$ and  $D^s(\gamma)$ bounded by $w_u(\gamma)$ and $w_s(\gamma)$, are equal:
$$
\int\limits_{D^u(\gamma)}\omega = \int\limits_{D^s(\gamma)}\omega.
$$
\end{Lm}
{\bf Proof}. The existence of stable manifold $W^s(\gamma)$ and unstable manifold $W^u(\gamma)$ of invariant
KAM-curves $\gamma \in W^c$ will be proved in Appendix.

The transversality condition implies the intersection of $W^{cu}$ with $W^{cs}$ near $q$ to occur along
a smooth 2-dimensional disk $\Sigma$. For every invariant curve $\gamma$ in $W^c$ its stable manifold
being extended by $f^{-1}$ in a finite number of iterations reaches a neighborhood of $q$ and transversely
intersects $\Sigma$ within $W^{cs}(\gamma)$ along closed curve $w_s(\gamma)$, the trace of $W^s(p)$ is point
$q$ itself. Traces on $\Sigma$ of $W^u(p)$ and $W^u(\gamma)$ in $W^{cu}(p)$ are respectively point $q$ and curve $w_u(\gamma)$.

Consider now a piece-wise smooth 2-dimensional surface made up of a piece of the lateral side of the cylinder $W^s(\gamma)$ between $W^c$ and $\Sigma$,
the piece of $W^c$ bounded by $\gamma$ and disk $\Sigma$. Integration of the form $\Omega$ over this surface is reduced to
the difference of integrals over the disk in $W^c$ and that over
disk in $\Sigma$ bounded by $w_s(\gamma)$, since the integral over lateral side is equal zero (it is
a Lagrangian submanifold). This gives the equality of integrals in the statement of the Lemma.
Similarly, we get equality of integral over disk in $W^c$, bounded by $\gamma$, and integral over disk in $\Sigma$,
bounded by $w_u(\gamma)$. $\blacksquare$

\section{Linearization and scattering map}\label{scattering_map_section}

The genericity Assumption \ref{scattering_as} is formulated using scattering map $S$. In this section
we will construct this map which acts on tangent plane $T_p W^c$. Scattering map is an analog of
the scattering matrix for a Schrodinger type equation \cite{Zakharov}. For the problems of the homoclinic dynamics
related with non-hyperbolic equilibria this map was first introduced in
\cite{KL1}. Far-reaching extension of this map for a normally hyperbolic
manifold in a Hamiltonian system was obtained in \cite{DLS}.

Consider first the linearization of the family of diffeomorphisms $f^n$ at homoclinic orbit
$\Gamma = \left\{q_n,\; n \in \N \; \big| \; q_{n+1} = f(q_n),\; q_0 = q \right\}$. This linearization is a sequence
of linear symplectic maps $L_n=Df \big|_{T_{q_n}M}: T_{q_n}M \rightarrow T_{q_{n+1}}M$ and hence $\lim L_n = Df_p$
as $|n|\to \infty$. Since $f^n(q)\to p$ as $n\to \pm\infty$, there exists an integer $N$ large enough such that
given a neighborhood $V$ of $p$ one gets $f^n(q) \in V$ for all $|n|\ge N$.

In neighborhood $V$ we choose a symplectic chart where fixed point $p$ is the
origin, then map $f$ is in the standard form "linear diffeomorphism plus higher order terms". After a linear
symplectic change of variables the linear part of the map can be transformed to the block-diagonal form:
\begin{equation}\label{linearized_map}
\left\{
\begin{array}{ccl}
x_{1} & = & \mu x + \ldots \\
y_{1} & = & \mu^{-1} y + \ldots \\
\left(
  \begin{array}{c}
    u_{1} \\ v_{1}
  \end{array}
\right) & = &
\left(
  \begin{array}{lr}
    \cos \alpha & -\sin \alpha \\
    \sin \alpha & \cos \alpha
  \end{array}
\right)
\left(
  \begin{array}{c}
    u \\ v
  \end{array}
\right) + \ldots,
\end{array}
\right.
\end{equation}
with $0< \mu <1$, dots mean terms of the order 2 and higher. In these coordinates the linearization
of this discrete dynamical system at the homoclinic orbit $\Gamma$ is given as follows:
\begin{equation} \label{linearized_system}
\left\{
\begin{array}{ccl}
\xi_{n+1} & = & \mu \xi_n + P_n \zeta_n \\
\eta_{n+1} & = & \mu^{-1} \eta_n + Q_n \zeta_n \\
\chi_{n+1} & = & R_{\alpha} \chi_n + W_n \zeta_n,
\end{array}
\right.
\end{equation}
where $\zeta_n = (\xi_n,\eta_n,\chi_n)^\top = (\xi_n,\eta_n,\chi^1_n,\chi^2_n)^\top$ is  coordinate
4-column vector in the tangent space at the point $q_n = (x_n,y_n,u_n,v_n)$; $R_{\alpha}$ denotes
the rotation matrix through angle $\alpha$:
$$
R_{\alpha} = \left(
  \begin{array}{lr}
    \cos \alpha & -\sin \alpha \\
    \sin \alpha & \cos \alpha
  \end{array}
\right);
$$
$P_n$, $Q_n$ are 1-row matrices, $W_n$ is $(4\times 2)$-matrix. Since $P_n, Q_n, W_n$ are of at least order 1
at $(0,0,0,0)$ and $(x_n,y_n,u_n,v_n)$ decay exponentially fast to $(0,0,0,0)$ as $|n|\to \infty$, for these matrices the following estimates
hold for $|n|\ge N$ and some positive $C$ depending on $N$ and on size of the neighborhood $V$:
$$
\left\|P_n\right\| \leq C \mu_1^{|n|}, \quad
\left\|Q_n\right\| \leq C \mu_1^{|n|}, \quad
\left\|W_n\right\| \leq C \mu_1^{|n|},
$$
where $0 < \mu < \mu_1 < 1$. Take $\left| n \right| \geq N$ and denote
$$
S_n = \left(
  \begin{array}{cc}
    E_2 & 0 \\
    0 & R_{n\alpha}
  \end{array}
\right),
$$
where $E_2$ is $2\times 2$ identity matrix.

Consider now the case $n\ge N$ and perform in the system (\ref{linearized_system}) a sequence of nonautonomous (with "time" $n$)
symplectic changes of variables $\zeta_n = S_{n-N}\psi_n$, where $\psi_n= (\hat{\xi}_n,\hat{\eta}_n,\hat{\chi}_n)$,
and consider (\ref{linearized_system}) in the rotating coordinate frame. This change of variables allows one to
exclude asymptotically the rotation in coordinates $\chi_n = (\chi^1_n,\chi^2_n)$ and prove that in new coordinates
each invariant bounded sequence for the linear system obtained from (\ref{linearized_system}) has a limit as $n \rightarrow \infty$.

After the change system (\ref{linearized_system}) casts as follows (we hold previous notations for variables):
\begin{equation} \label{linearized_wo_rotation}
\left\{
\begin{array}{ccl}
\xi_{n+1} & = & \mu \xi_n + F_n \zeta_n \\
\eta_{n+1} & = & \mu^{-1} \eta_n + G_n \zeta_n \\
\chi_{n+1} & = & \chi_n + H_n \zeta_n,
\end{array}
\right.
\end{equation}
where $F_n,\; G_n$ are again 1-row matrices and $H_n$ is $4\times 2$-matrix. For these matrices estimates
similar to those for matrices $P_n$, $Q_n$ and $W_n$ are valid. Sequence $\left\{(\xi_n,\eta_n,\chi_n)\right\}$
is called the solution of the system or
the invariant sequence, if equalities (\ref{linearized_wo_rotation}) are satisfied for all $n\in \Z_+$. The following lemma is valid.
\begin{Lm}\label{integral_system_solution_lm}
There is an integer $N > 0$ large enough such that for any given $\xi^0\in \R$, $\chi_{+} \in \R^2$
a unique solution $\left\{(\xi_n,\eta_n,\chi_n)\right\}$, $n \geq N$, exists for the system $(\ref{linearized_wo_rotation})$
such that for this solution the boundary conditions are satisfied: $\xi_N = \xi^0$, $\chi_n \rightarrow \chi_{+}$, $\left|\xi_n\right| \rightarrow 0$,
$\left|\eta_n\right| \rightarrow 0$ as $n \rightarrow + \infty$.
\end{Lm}
{\bf Proof}. Similar to \cite{KL2}, instead of system (\ref{linearized_wo_rotation}) consider a system of
difference equations (\ref{integral_system}):
\begin{equation}\label{integral_system}
\left\{
\begin{array}{ccl}
\xi_n & = & \mu^{n-N}\xi^0 + \sum\limits_{s=N}^{n-1} \mu^{n-1-s}F_s\zeta_s \\
\eta_n & = & -\sum\limits_{s=n}^{\infty} \mu^{s+1-n}G_s\zeta_s \\
\chi_n & = & \chi_{+} - \sum\limits_{s=n}^{\infty} H_s\zeta_s.
\end{array}
\right.
\end{equation}
Note that any solution of this system obeys the boundary conditions in the statement of the lemma. Let us show
first that the solution of the system (\ref{integral_system}) is also the solution of the system
(\ref{linearized_wo_rotation}) and vice versa. Indeed, the following equalities hold:
\begin{multline*}
\xi_{n+1} = \mu^{n+1-N} \xi^{0} + \sum\limits^{n}_{s=N} \mu^{n-s} F_s \zeta_s \\[\jot]
\shoveright{= \mu\left(\mu^{n-N} \xi^{0} + \sum\limits^{n-1}_{s=N} \mu^{n-1-s}F_s \zeta_s\right) + F_{n}\zeta_n
            = \mu \xi_n + F_{n}\zeta_n,}\\
\shoveleft{\eta_{n+1} = -\sum\limits^{\infty}_{s=n+1} \mu^{s-n}G_s \zeta_s
           =  G_n\zeta_n - G_n\zeta_n -\sum\limits^{\infty}_{s=n+1} \mu^{s-n}G_s\zeta_s}\\
\shoveright{= G_n\zeta_n + \mu^{-1}\left( -\sum\limits^{\infty}_{s=n}
            \mu^{s+1-n}G_s\zeta_s\right)=\mu^{-1}\eta_n+G_n\zeta_n,}\\
\shoveleft{\chi_{n+1} = \chi_{+} - \sum\limits_{s=n+1}^{\infty} H_s \zeta_s
           = \chi_{+} - \sum\limits_{s=n}^{\infty} H_s \zeta_s + H_n \zeta_n
           = \chi_n + H_n \zeta_n.}
\end{multline*}
So, if the sequence $\left\{(\xi_n,\eta_n,\chi_n)\right\}$ solves (\ref{integral_system}), then it satisfies
(\ref{linearized_wo_rotation}). The converse assertion is given as $n\to \infty$ by the consecutive application of
(\ref{linearized_wo_rotation}) to an initial point.

Thus, one needs to prove the existence of solutions for system (\ref{integral_system}). To do this, we use
the contraction mapping principle. Denote $B$ the Banach space of sequences $\zeta = \left\{(\xi_n,\eta_n,\chi_n)\right\}$
uniformly bounded on $[N, + \infty)$ with the norm
$$
\left\| \zeta \right\| = \underset{n\geq N}{\sup}(\left|\xi_n\right|, \left|\eta_n\right|, \left\|\chi_n\right\|).
$$
Right hand sides of (\ref{integral_system}) define operator $T:\; T[\zeta] = \bar{\zeta}$ on $B$. At the first step
let us verify that $T$ is defined correctly, that is $T[B] \subset B$, here $\chi_{+}$, $\xi^{0}$ and $N$ are considered
as parameters. Recall that for $F_n$, $G_n$, $H_n$ the following estimates are valid:
$\left\|F_n\right\|, \left\|G_n\right\|, \left\|H_n\right\| \leq C \nu^n$, $0<\nu<1$. Here $C$ depends on $N$, but is
finite for a fixed $N$. Denote $\kappa = \max \left\{\mu, \nu \right\} <1$. Then one proceeds as follows:
\begin{multline*}
\left| \bar{\xi}_{n} \right|
     \leq \kappa^{n-N} \xi^{0} + \left\| \zeta \right\| \sum\limits^{n-1}_{s=N} \kappa^{n-1-s} C \kappa^{s}
     = \kappa^{n-N} \xi^{0} + \left\| \zeta \right\| C \kappa^{n-1}(n-N),\\
\shoveleft{\left| \bar{\eta}_{n} \right|
     \leq \left\| \zeta \right\| \sum\limits^{\infty}_{s=n} \kappa^{s+1-n} C \kappa^s
     = C \left\| \zeta \right\| \kappa^{1-n} \sum\limits^{\infty}_{s=n} \kappa^{2s}
     = C \left\| \zeta \right\| \kappa^{1-n} \displaystyle{\frac{\kappa^{2n}}{1-\kappa^2}}}\\
\shoveright{\displaystyle{ = C \left\|\zeta \right\| \frac{\kappa^{n+1}}{1-\kappa^2}},}\\
\shoveleft{ \left\| \bar{\chi}_{n} \right\|
     \leq \left\| \chi_{+} \right\| + C \left\| \zeta \right\| \sum\limits_{s=n}^{\infty} \kappa^{s}
     = \left\| \chi_{+} \right\| + \displaystyle{C \left\| \zeta \right\| \frac{\kappa^n}{1-\kappa}}.}\\
\end{multline*}
Thus, the sequence $\bar{\zeta}$ is uniformly bounded on $[N, + \infty)$, so the operator $T$ is defined correctly.

Next we prove $T$ to be a contraction map:
\begin{multline*}
\left| \bar{\xi}_{n}^1 - \bar{\xi}_{n}^2 \right|
           \leq \left\| \zeta_1 - \zeta_2 \right\| C \sum\limits^{n-1}_{s=N} \kappa^{n-1-s}\kappa^s
           = \left\| \zeta_1 - \zeta_2 \right\| C \kappa^{n-1}(n-N),\\
\shoveleft{\left| \bar{\eta}_{n}^1 - \bar{\eta}_{n}^2 \right|
           \leq \left\| \zeta_1 - \zeta_2 \right\| C \sum\limits^{\infty}_{s=n} \kappa^{s+1-n} \kappa^s
           = \left\| \zeta_1 - \zeta_2 \right\| C \kappa^{1-n} \displaystyle{\frac{\kappa^{2n}}{1-\kappa^2}}}\\
\shoveright{= \left\| \zeta_1 - \zeta_2 \right\| \displaystyle{\frac{C\kappa^{n+1}}{1-\kappa^2}}
            \leq \left\| \zeta_1 - \zeta_2 \right\| \displaystyle{\frac{C\kappa^{N}}{1-\kappa^2}},}\\
\shoveleft{\left\| \bar{\chi}_{n}^1 - \bar{\chi}_{n}^2 \right\|
           \leq \left\| \zeta_1 - \zeta_2 \right\| C \sum\limits_{s=n}^{\infty} \kappa^{s}
           = \left\| \zeta_1 - \zeta_2 \right\| \displaystyle{\frac{C\kappa^n}{1-\kappa}}
           \leq \left\| \zeta_1 - \zeta_2 \right\| \displaystyle{\frac{C\kappa^N}{1-\kappa}}.}
\end{multline*}
These estimates show that $T$ is contracting for $N$ large enough and $n \geq N$. Thus, for any fixed $\xi^0 \in \R$
and  $\chi_{+} \in \R^2$ there is a unique solution $\zeta(\xi^0, \chi_{+}) = \{\zeta_n(\xi^0, \chi_{+})\} =
\{(\xi_n(\xi^0, \chi_{+}), \eta_n(\xi^0, \chi_{+}), \chi_n(\xi^0, \chi_{+}))\}$ for the system (\ref{integral_system})
such that $\xi_N(\xi^0, \chi_{+}) = \xi^0$. The estimates above also show that $\left| \xi_n \right|$,
$\left| \eta_n \right|$ and $\left\| \chi_n - \chi_{+} \right\|$ tend to zero as $n \rightarrow \infty$.
$\blacksquare$

For the further purposes one needs to prove some linearity relations for solutions of system (\ref{linearized_wo_rotation}).
\begin{Lm}\label{linearity_relations}
Solutions $\zeta(\xi^0,\chi_{+})$ of the system $(\ref{linearized_wo_rotation})$ satisfy the following linearity relations:
\begin{enumerate}[I.] 
\item $\zeta(\xi'_0,\chi_{+}) + \zeta(\xi^{\prime\prime}_0,0) =
       \zeta(\xi'_0 + \xi^{\prime\prime}_0,\chi_{+})$;
\item $\zeta(\xi'_0,0) + \zeta(\xi^{\prime\prime}_0,0) =
       \zeta(\xi'_0 + \xi^{\prime\prime}_0,0)$, \quad
      $\zeta(\alpha \xi'_0,0) = \alpha \zeta(\xi'_0,0)$;
\item $\zeta(0,\alpha \chi_{+} + \beta \chi^{\prime\prime}_{+}) =
       \alpha \zeta(0,\chi'_{+}) + \beta \zeta(0,\chi^{\prime\prime}_{+})$.
\end{enumerate}
\end{Lm}
{\bf Proof}. To prove the first equality consider the function
\begin{multline*}
\Delta_1 = \left\{ \Delta_1(n) \right\} =
           \zeta(\xi'_0,\chi_{+}) +
           \zeta(\xi^{\prime\prime}_0,0) -
           \zeta(\xi'_0 + \xi^{\prime\prime}_0,\chi_{+}) \\
= \left\{ \left(
  \begin{array}{c}
			\xi_n(\xi'_0,\chi_{+}) +
      \xi_n(\xi^{\prime\prime}_0,0) -
      \xi_n(\xi'_0 + \xi^{\prime\prime}_0,\chi_{+})\\
			\eta_n(\xi'_0,\chi_{+}) +
      \eta_n(\xi^{\prime\prime}_0,0) -
      \eta_n(\xi'_0 + \xi^{\prime\prime}_0,\chi_{+})\\
      \chi_n(\xi'_0,\chi_{+}) +
      \chi_n(\xi^{\prime\prime}_0,0) -
      \chi_n(\xi'_0 + \xi^{\prime\prime}_0,\chi_{+})
	\end{array}
\right) \right\}.
\end{multline*}
This function is a solution of the system (\ref{integral_system}) with boundary conditions $(0,0)$. Indeed,
consider the following systems with boundary conditions $(\xi'_0,\chi_{+})$, $(\xi^{\prime\prime}_0,0)$ and
$(\xi'_0 + \xi^{\prime\prime}_0,\chi_{+})$, respectively:
\begin{equation}\label{system_1}
\left\{
\begin{array}{ccl}
\xi_n(\xi'_0,\chi_{+}) & = &
        \mu^{n-N}\xi^0 + \sum\limits_{s=N}^{n-1} \mu^{n-1-s}F_s\zeta_s(\xi'_0,\chi_{+}) \\
\eta_n(\xi'_0,\chi_{+}) & = &
        -\sum\limits_{s=n}^{+\infty} \mu^{s+1-n}G_s\zeta_s(\xi'_0,\chi_{+}) \\
\chi_n(\xi'_0,\chi_{+}) & = &
        \chi_{+} - \sum\limits_{s=n}^{+\infty} H_s\zeta_s(\xi'_0,\chi_{+}),
\end{array}
\right.
\end{equation}

\begin{equation}\label{system_2}
\left\{
\begin{array}{ccl}
\xi_n(\xi^{\prime\prime}_0,0) & = &
      \mu^{n-N}\xi^0 + \sum\limits_{s=N}^{n-1} \mu^{n-1-s}F_s\zeta_s(\xi^{\prime\prime}_0,0) \\
\eta_n(\xi^{\prime\prime}_0,0) & = &
      -\sum\limits_{s=n}^{+\infty} \mu^{s+1-n}G_s\zeta_s(\xi^{\prime\prime}_0,0) \\
\chi_n(\xi^{\prime\prime}_0,0) & = &
      \chi_{+} - \sum\limits_{s=n}^{+\infty} H_s\zeta_s(\xi^{\prime\prime}_0,0),
\end{array}
\right.
\end{equation}

\begin{equation}\label{system_3}
\left\{
\begin{array}{ccl}
\xi_n(\xi'_0 + \xi^{\prime\prime}_0,\chi_{+}) & = &
     \mu^{n-N}\xi^0 + \sum\limits_{s=N}^{n-1}
     \mu^{n-1-s}F_s\zeta_s(\xi'_0 + \xi^{\prime\prime}_0,\chi_{+}) \\
\eta_n(\xi'_0 + \xi^{\prime\prime}_0,\chi_{+}) & = &
     -\sum\limits_{s=n}^{+\infty}
     \mu^{s+1-n}G_s\zeta_s(\xi'_0 + \xi^{\prime\prime}_0,\chi_{+}) \\
\chi_n(\xi'_0 + \xi^{\prime\prime}_0,\chi_{+}) & = &
     \chi_{+} - \sum\limits_{s=n}^{+\infty}
     H_s\zeta_s(\xi'_0 + \xi^{\prime\prime}_0,\chi_{+}).
\end{array}
\right.
\end{equation}

Summing related equalities from (\ref{system_1}), (\ref{system_2}) and subtracting (\ref{system_3}) we get:
\begin{multline*}
\xi_n(\xi'_0,\chi_{+}) + \xi_n(\xi^{\prime\prime}_0,0)
    - \xi_n(\xi'_0 + \xi^{\prime\prime}_0,\chi_{+})
    = \mu^{n-N}(\xi'_0 + \xi^{\prime\prime}_0
    - (\xi'_0 + \xi^{\prime\prime}_0)) \\
+ \sum\limits^{n-1}_{s=N} \mu^{n-s-1} F_s(\zeta_s(\xi'_0,\chi_{+})
    + \zeta_s(\xi^{\prime\prime}_0,0)
    - \zeta_s(\xi'_0 + \xi^{\prime\prime}_0,\chi_{+}))\\
= \sum\limits^{n-1}_{s=N} \mu^{n-s-1} F_s \Delta_1(s),
\end{multline*}
\begin{multline*}
\eta_n(\xi'_0,\chi_{+}) + \eta_n(\xi^{\prime\prime}_0,0)
    - \eta_n(\xi'_0 + \xi^{\prime\prime}_0,\chi_{+})\\
= -\sum\limits^{+\infty}_{s=n} \mu^{-(n-s-1)} G_s(\zeta_s(\xi'_0,\chi_{+})
    + \zeta_s(\xi^{\prime\prime}_0,0)
    - \zeta_s(\xi'_0+ \xi^{\prime\prime}_0,\chi_{+})) \\
= -\sum\limits^{+\infty}_{s=n} \mu^{-(n-s-1)} G_s \Delta_1(s),
\end{multline*}
\begin{multline*}
\chi_n(\xi'_0,\chi_{+}) + \chi_n(\xi^{\prime\prime}_0,0)
    - \chi_n(\xi'_0 + \xi^{\prime\prime}_0,\chi_{+}) \\
= (\chi_{+} + 0 - \chi_{+}) -	\sum\limits_{s=n}^{+\infty}
    H_s (\zeta_s(\xi'_0,\chi_{+})
    + \zeta_s(\xi^{\prime\prime}_0,0)
    - \zeta_s(\xi'_0 + \xi^{\prime\prime}_0,\chi_{+}))\\
= -	\sum\limits_{s=n}^{+\infty} H_s \Delta_1(s).
\end{multline*}
The equalities obtained imply that $\Delta_1$ satisfies (\ref{integral_system}) with boundary conditions
$(\xi^0,\chi_{+})=(0,0)$. Since the solution of (\ref{integral_system}) with given boundary conditions is
unique, then $\Delta_1 \equiv 0$ and therefore the relation $I$ is valid.

Relations $II$ and $III$ are proved in a similar way, if instead of $\Delta_1$ one considers $\Delta_2$
and $\Delta_3$, respectively:
$$
\Delta_2 =
       \zeta(\alpha \xi'_0 + \beta \xi^{\prime\prime}_0, 0) -
       \alpha \zeta(\xi'_0, 0) -
       \beta \zeta(\xi^{\prime\prime}_0, 0),
$$
$$
\Delta_3 =
      \zeta(0, \alpha \chi'_{+} + \beta \chi^{\prime\prime}_{+}) -
      \alpha \zeta(0,\chi'_{+}) -
      \beta \zeta(0,\chi^{\prime\prime}_{+}).
$$
Lemma has been proved. $\blacksquare$

Similar lemmas hold for $n \leq -N$.

\subsection{Geometry of linearized map}

Now let us present a geometrical interpretation of the results obtained.
To do this, we introduce a countable set of linear symplectic spaces
$X_n,$ $|n|\ge N,$ with coordinates $(\xi_n,\eta_n,\chi^1_n,\chi^2_n),$ and linear symplectic maps
$\mathcal L_n: X_n \to X_{n+1},$ $n\ge N,$ defined by
(\ref{integral_system}).
If we fix $\chi_+$ and vary $\xi_0$, $\xi_0 = \xi_0^{\prime} + \xi_0^{\prime\prime},$ $\xi_0^{\prime\prime} \in \R,$ then
due to Lemma \ref{linearity_relations}, the related solutions of (\ref{integral_system}) define an affine straight
line in $X_N$ (in fact, they are initial points of these solutions) and hence in any $X_n,$ $n> N.$ These straight lines have
the characteristic property that any solution which passes through this
line in $X_n$ decays exponentially as $n\to \infty$: $\xi_n \to 0,$ $\eta_n \to 0,$ $||\chi_n - \chi_+|| \to 0.$ In addition, if we fix not $\chi_+$
but only the value
$2I = \chi_1^2 + \chi_2^2 = ||\chi_+||^2$, then in every $X_n,$ $n\ge N$ we get a 2-dimensional
cylinder $C_n^+(I)$ formed of those straight lines in $X_n$ through which
solutions asymptotically satisfy to $||\chi_n - \chi_+|| \to 0,$ each
straight line on $C_n^+(I)$ corresponds to the unique value on the circle $2I = \chi_1^2 +
\chi_2^2$ (an asymptotic phase). Varying $I, \xi_0$ defines a linear 3-dimensional subspace $\mathcal L_n^+$ of bounded solutions in $X_N$,
and hence in $X_n$, which in turn foliates on the cylinders $C_n^+(I)$. Such a cylinder shrinks
to the straight line $C_n^+(0)$, as $I \to 0$, this straight line just corresponds to solutions with $\chi_+ =0$.

Now let us turn to the initial linearization problem along homoclinic orbit $\Gamma$ for diffeomorphism
$f$. To derive results described above, we performed the sequence of linear changes of
variables that allowed us to prove for any bounded solution the existence
of an asymptotic phase. In the initial coordinates all objects found
preserve: for any point $q_n$ in the related tangent space $T_{q_n}M$ we
have 3-dimensional subspace of bounded as $n\to \infty$ solutions $\mathcal
L_n^+$ (we preserve the same notations for similar objects) which are
foliated into cylinders $C_n^+(I)$ (it is worth mentioning that the value of $I$
does not change when returning to the initial coordinates), foliations
into straight lines, etc. It is evident that in fact $\mathcal
L_n^+$ is nothing else as tangent space $T_{q_n}W^{cs}$

The same picture takes place for $X_n$ with $n\le -N$, the only difference
is that one needs take limits as $n\to -\infty.$ Here we also have
cylinders $C_n^-$, straight lines, 3-subspaces $\mathcal L_n^- = T_{q_n}W^{cu}$, and so forth.

For tangent space $T_{q(-N)}M$ and $T_{q(N)}M$ we have linear symplectic map $\mathcal S =
Df^{2N}: T_{q(-N)}M \to T_{q(N)}M$ calculated at the point $q(-N).$ This map transforms $\mathcal
L_n^-$ to a 3-dimensional subspace in $T_{q(N)}M$ which transversely intersects the
straight line $C_N^+(0)= l_N^s$ being the tangent space to $W^s(p)$.

\subsection{Scattering map}

\begin{figure}
\centering
\parindent=0pt
\includegraphics[width=0.9\textwidth]{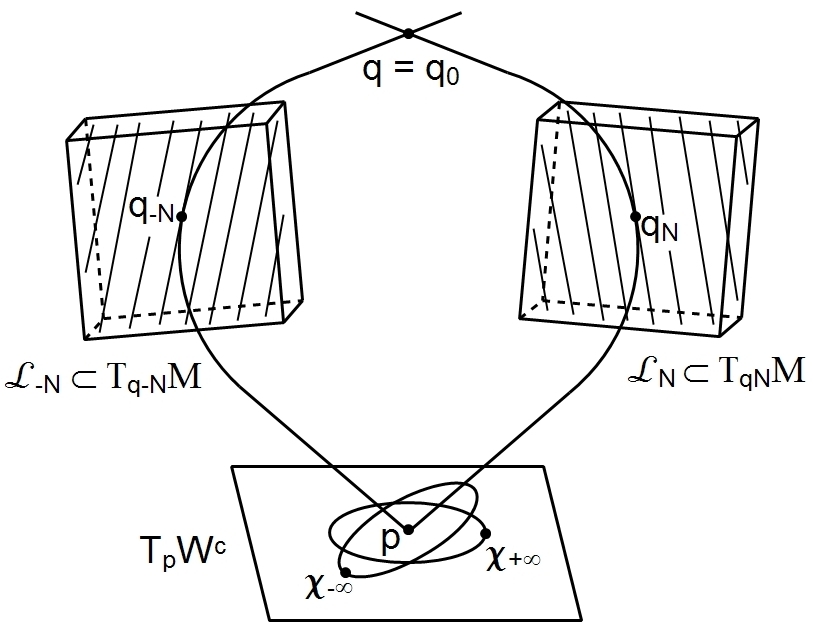}
\caption{Scattering map construction}
\label{fig:scattering_map}
\end{figure}

Now we are ready to construct the scattering map $S: T_p W^c \rightarrow T_p W^c$. Take any point $s \in T_p W^c$.
Fixing this point defines the unique straight line in $T_{q(-N)}M$ of the foliation defined in $\mathcal{L}_{-N}=T_{q(-N)}W^{cu}$ whose points are
asymptotic to $s$ as $n \to -\infty$. Let us apply linear map $Df^{2N}$ to points of this line. We get the straight
line in $T_{q(N)}M$ which is transversal to 3-plane $\mathcal{L}_N=T_{q(N)}W^{cs}$ due to transversality condition. Thus, the line
obtained intersects this 3-plane at the unique point through which a unique line of the foliation defined in plane
$\mathcal{L}_N$ passes. Denote $s_1 \in T_p W^c$ that unique point which is the limit as $n\to \infty$ for all
sequences starting on this line. We set $S(s) = s_1$ (Fig.~\ref{fig:scattering_map}).

Let us verify that $S$ is a linear map. It is clear that $S(0)=0$. Indeed, for $s = 0$ the corresponding
straight line in 3-plane $\mathcal L_{-N}$ is the tangent line to $W^u$ in tangent space $T_{q(-N)}M$.
Its image under $\mathcal S =Df^{2N}$ is a straight line in $T_{q(N)}M$ which is transversal to $T_{q(N)}W^{cs}$ due
to Assumption \ref{transversality_as} and intersects it at the origin of $T_{q(N)}M$. Through the origin
the unique line from the constructed foliation passes: the tangent line to $W^s$ which corresponds to $s_1 = 0$ in $T_pW^c$.

Denote $l_{s}$ that straight line in $\mathcal{L}_{-N}$ which consists of points through which
solutions pass tending to $s$ as $n \rightarrow -\infty$. Using the linearity relations \textit{I} and \textit{II}
we get the following representation for the solutions of system (\ref{linearized_wo_rotation}):
$$
\zeta_n(\xi^0,\lambda \chi_{-}) = \zeta_n(0,\lambda \chi_{-}) + \zeta_n(\xi^0,0) =
\lambda \zeta_n(0,\chi_{-}) + \zeta_n(\xi^0,0)\;\mbox{\rm for any}\; \lambda \in \mathbb R.
$$
To find the image $S(\lambda s)$ we act by $Df^{2N}$ on $\lambda l_{s} + l_0^u$ (here $l_0^u$
is the tangent line to $W^u$ at point $q(-N)$). In $\mathcal{L}_N$ we get two vectors, through each such
vector a unique straight line passes: these lines are $\lambda l_{s}$ and $l_0^s$, respectively ($l_0^s$ is
tangent line to $W^s$ at point $q_N$). Thus, $S$ maps $\lambda s$ to $\lambda s_1 + 0 = \lambda s_1$.
Similarly, for the sum $s' + s^{\prime\prime}$ we get relations for corresponding $\chi'_{-}, \chi^{\prime\prime}_{-}$:
\begin{multline*}
\zeta_n(\xi^0,\chi'_{-} + \chi^{\prime\prime}_{-}) =
\zeta_n(\xi^0 + 0,\chi'_{-} + \chi^{\prime\prime}_{-}) =\\
=\zeta_n(\xi^0,0) + \zeta_n(0,\chi'_{-} + \chi^{\prime\prime}_{-}) =
\zeta_n(\xi^0,0) + \zeta_n(0,\chi'_{-}) + \zeta_n(0,\chi^{\prime\prime}_{-}).
\end{multline*}
In this case we act by $Df^{2N}$ on $l_{s'} + l_{s^{\prime\prime}} + l_0^u$.

The next proposition characterizes map $S$.
\begin{Prop}
Map $S: T_pW^c \to T_pW^c$ is symplectic.
\end{Prop}
{\bf Proof}.
Choose any two vectors in  the symplectic plane $T_pW^c$. These vectors
define two straight lines from the foliation in $\mathcal L_{-N}$. Take then two vectors $v_1, v_2$ in $\mathcal L_{-N} = T_{q(-N)}W^{cu}$
corresponding to these lines: origins of vectors coincide with zero point of $T_{q(-N)}W^{cu}$ and ends of the vectors
belong to corresponding line. Skew-scalar product being the restriction of 2-form  $\Omega$ on tangent space $\mathcal L_{-N}$,
does not depend on vectors we choose. Indeed, difference of vectors corresponding to the same line is vector lying in
$l^u_{-}$ which is zero vector for skew-scalar product (such vectors shrink exponentially in backward iterations).
$\mathcal S$-images of these two straight lines are two straight lines in $T_{q(N)}M$ which are transversal to subspace $\mathcal L_N$.
The intersection of the lines with $\mathcal L_N$ gives two vectors $T(v_1), T(v_2)$, whose origins coincides with zero of
$\mathcal L_N = T_{q(N)}W^{cs}$, for specified $v_1, v_2$. Since $\mathcal S$ is linear symplectic map, then the skew-scalar product is preserved.
Now we have two straight lines from foliation in $mathcal L_N$ and again skew-scalar product of vectors corresponding to different
lines does not depend on exact vectors we choose. But this product is equal to skew-scalar product of vectors $v_1, v_2$
and does not change in forward iterations. Therefore skew-scalar product in limit in $T_pW^c$ equals to skew-scalar
product of initial vectors in $T_pW^c$. $\blacksquare$

Linear symplectic map $S$ we call \textit{the scattering map}.

\section{Homoclinic orbits to invariant KAM-curves}\label{homoclinic_orbits_existence_proof}

\begin{figure}
\centering
\parindent=0pt
\includegraphics[width=0.9\textwidth]{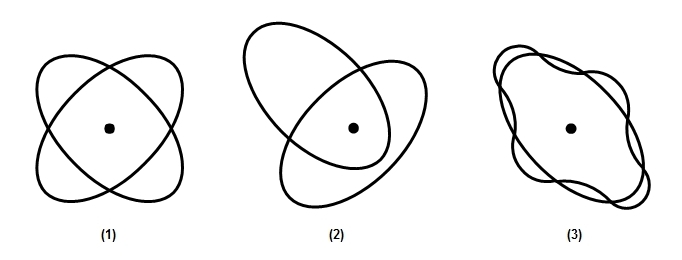}
\caption{Possible intersections of $W^s(\gamma)$, $W^u(\gamma)$ on $\Sigma$}
\label{fig:intersections}
\end{figure}

To prove Theorem \ref{main_th} we assume for diffeomorphism $f$ Assumptions \ref{intersection_as}, \ref{transversality_as}
and  \ref{scattering_as} to hold. Thus, according to Section \ref{transversality_condition_consequence} manifolds
$W^{cu}(p)$ and $W^{cs}(p)$ intersect at homoclinic point $q$ transversally, and therefore along a symplectic 2-disk $\Sigma$
(we may regard $q=q_+$ and disk $\Sigma$ to belong to a neighborhood of fixed point $p$). The idea of the existence proof
for homoclinic orbits to invariant KAM-curve $\gamma \in W^c$ is the following. Let $V$ be a sufficiently
small neighborhood of point $p$. For each KAM-curve $\gamma \subset W^c\cap V$ its action is defined according
to the Stokes theorem as the integral of 2-form $\Omega$ over that disk in $W^c$ whose boundary is curve $\gamma$.
Curve $\gamma$ has its local stable manifold which can be extended by a finite number of iterations of map $f^{-1}$ till
the manifold reaches the neighborhood of homoclinic point $q$ staying inside of $W^{cs}(p)$. Therefore, this manifold
intersects transversally within $W^{cs}(p)$ disk $\Sigma$ along a closed curve $w_s(\gamma)$. Similarly, unstable manifold of the same
curve $\gamma$ under the action of $f$ reaches the neighborhood of $q$ staying inside $W^{cu}(p)$ and hence intersects $\Sigma$
along a closed curve $w_u(\gamma)$. Two obtained curves on $\Sigma$ have the same value of action as follows from lemma
\ref{equal}. Thus, two disks in $\Sigma$ bounded by $w_s(\gamma),$ $w_u(\gamma)$ are of the same area and have common
point $q$ lying inside both of them.
Hence, the intersection of curves $w_s(\gamma)$ and $w_u(\gamma)$ is not empty and consists of at least two different points,
homoclinic orbits to $\gamma$ pass through the intersection points (see Fig.~\ref{fig:intersections}).

The problem here is that we do not know a precise information on this intersection: how many
points does it contain, if the intersection is transverse or not, etc. All these questions are relevant
for the further study of nearby dynamics. In the case of an integrable diffeomorphism
(when an invariant w.r.t. the diffeomorphism $f$ smooth function exists)
these two curves coincide, since both of them belong to the same level of the invariant function whose restriction
on $\Sigma$ near $q$ usually forms a connected closed curve.

Provided that our Assumptions hold, we shall prove for
any KAM curve on $W^c$ in a sufficiently small $V$ with a given value of action $I$, the intersection to consist
of exactly four points and it is transverse at each of these points. To prove this we shall connect intersection
properties for curves $w_u(\gamma)$ and $w_s(\gamma)$ on $\Sigma$ with intersection properties of related ellipses
$E$ and $S(E)$ with the same action $I$ in tangent plane $T_p W^c$.
The genericity Assumption implies that these ellipses have the same center and the same area and intersect
transversally at exactly four points. This property will be carried to the intersection of $w_u(\gamma)$ and $w_s(\gamma)$.

To prove the intersection of theses curve as in Fig.~\ref{fig:intersections} (left panel), we transform $f$ in neighborhood $V$
of fixed point $p$ to the normal form (\ref{normal_diffeomorphism}) (see Appendix) up to the third order terms and consider first the truncated map.
For this map two local functions $xy$ and $u^2 + v^2$ are local integrals up to sixth order terms,
but what is more important for our goals, circles $x=y=0,$ $u^2 + v^2 = c$ on $W^c$ are invariant curves for any
positive $c$ small enough and stable and unstable manifolds of these curves have the representation $y=0$,
$u^2+v^2 = c$ and $x=0$, $u^2+v^2 = c$, respectively. Since the whole map differs from the truncated map
by terms of the fourth order and higher, then for a given positive constant $c$ small enough invariant manifolds of KAM-curves
for the truncated and full maps are at least $C^3$-close in $V$, due to Fenichel theorems. According to these theorems \cite{Fenichel2},
as diffeomorphism is $C^r$-smooth, then invariant cylinders of KAM-curves are $C^{r-1}$-smooth (see
Appendix). Also we shall suppose, without loss of generality, that both points
$q_- = q(-N)$ and $q_+ = q(N)$ belong to $V$.

For diffeomorphism (\ref{normal_diffeomorphism}) for $n\ge N$, due to properties of functions $\varphi,\,\psi$,
system (\ref{linearized_wo_rotation})  has the linearization matrices for $f$ along the homoclinic
orbit
$$
\begin{pmatrix}
\mu(1+\cdots)&p_n &q_n &r_n\\0&\mu^{-1}(1+\cdots)&0&0\\0&l_n &\cos\alpha +\cdots &-\sin\alpha +
\cdots\\0&s_n &\sin\alpha +\cdots &\cos\alpha +\cdots
\end{pmatrix},
$$
where dots and  $p_n, q_n, r_n, l_n, s_n$ mean terms tending to zero exponentially fast as $n\to +\infty$.
The form of this matrix implies that 3-dimensional plane $\eta =0$ in the tangent space at homoclinic point
$q(N)$ under the action of this linear map is transformed to 3-dimensional plane $\eta =0$
in the tangent space at homoclinic points $q(N+1)= f(q(N))$, etc. As Lemma
\ref{integral_system_solution_lm} implies, for a fixed $I$, $(\chi^1)^2 + (\chi^2)^2 = 2 I$,
we get in $T_{q(N)}M$ a cylinder in 3-dimensional plane $\eta =0$,
 consisting of solutions for the system (\ref{linearized_wo_rotation}) which asymptotically tend to the circle
$(\chi^1)^2 + (\chi^2)^2 = 2 I$ in $T_pW^c$. Intersection of this cylinder
with tangent plane to $\Sigma_+$ at point $q(N)$ is an ellipse. Similar cylinders and ellipses are
obtained, if one considers linearization of diffeomorpsim along homoclinic orbit for $n\le -N$ as $n\to -\infty$.

Global symplectic map $\mathcal{G}=f^{2N}$ transforms a neighborhood of homoclinic point $q(-N)$ to
a neighborhood of homoclinic point $q(N)$. In normalized coordinates (\ref{normal_diffeomorphism}), homoclinic points have coordinates:
$q(-N)=(0,y_1, 0, 0)$, $q(N)=(x_0, 0, 0, 0)$. Therefore, symplectic map $\mathcal{G}$ has the following
local representation
$$
\begin{array}{ccl}
\bar{x} - x_0 & = & F(x, y-y_1,u,v),\\
\bar{y} & = & G(x, y-y_1,u,v),\\
\bar{u} & = & P(x, y-y_1,u,v),\\
\bar{v} & = & Q(x, y-y_1,u,v),
\end{array}
$$
$$
d\bar{x}\wedge d\bar{y} + d\bar{u}\wedge d\bar{v} = dx\wedge dy + du\wedge dv,
$$
and all functions $F,G,P,Q$ vanish at the point $(0,y_1,0,0).$
Transversality condition at point $q(N)$ means the tangent vector to $W^{s}$, that is $
(1,0,0,0)$, be transverse to tangent plane to $W^{cu}$ being $\mathcal{G}$-image of 3-disk $x=0$.
This implies that determinant
$$
\displaystyle{ \det \left( \frac{\partial (G,P,Q)}{\partial (y-y_1,u,v)} \right)}
$$
calculated at point $(0,y_1,0,0)$ does not vanish. In virtue of Lemma \ref{conv} the same transversality condition
holds at $q(-N)$: tangent vector to $W^u$ (i.e. $(0,1,0,0)$) and tangent plane to $W^{cs}$ ($\mathcal{G}$-pre-image of $\bar{y}=0$)
are transverse, this is equivalent to inequality $G_y \ne 0$ at $(0,y_1,0,0)$.
2-disk $\Sigma = \Sigma_+$ at these coordinates is the intersection of $\mathcal{G}$-image of local 3-disk $x=0$ near point $q(-N)$ and
local 3-disk $\bar{y}=0$ near $q(N)$. Then one has a  representation for  $\Sigma_- = \mathcal{G}^{-1}(\Sigma_+)$ in the form
$x=0,\;y-y_1 = \Psi(u,v)$, and for $\Sigma_+$ in the form $\bar{x} = \Phi(\bar{u},\bar{v}),\;\bar{y}=0$.
In particular, since map $\mathcal{G}$ is symplectic and the restrictions of symplectic 2-form $\Omega$ to
symplectic disks $\Sigma_-$ and $\Sigma_+$ are $du\wedge dv$ and $d\bar{u}\wedge d\bar{v},$
respectively, then the restriction $\mathcal G_{\Sigma}$ of map $\mathcal{G}$, $\mathcal G_{\Sigma}: \Sigma_- \to \Sigma_+$
is symplectic and has the form
\begin{equation}\label{restr}
\begin{array}{l}
\bar{u}= a_{11}u + a_{12}v+ \ldots,
\bar{v}= a_{21}u + a_{22}v+ \ldots,
\end{array}
\end{equation}
where dots mean  terms of the second order and higher, matrix
$$
A=\begin{pmatrix}
a_{11}&a_{12}\\a_{21}&a_{22}
\end{pmatrix}
$$
is symplectic, here this means to be uni-modular: $\det A = 1$. In fact,
matrix $A$ depends on integer parameter $N$, $A = A_N$, since our choice of the
homoclinic points $q(-N),q(N)$ and hence $\Sigma_-,$ $\Sigma_+,$
and related tangent planes to them depends on $N$. But for any $N$ large
enough the representations for these tangent plane are similar and
coordinates of them are coordinates on $T_pW^c$. As $N\to \infty$ the
related tangent planes to $\Sigma_+(N)$ tend to $T_pW^c$. This implies
\begin{Lm}
For $N$ large enough matrix $A_N$ is not a rotation matrix.
\end{Lm}
{\bf Proof}. To prove Lemma we shall show that its conclusion follows from the transversality
Assumption \ref{transversality_as} and genericity Assumption
\ref{scattering_as}. As follows from Lemmas 3, 4, $A_N$ tends to matrix
$A_\infty$ being the coordinate representation of scattering map $S$. $\blacksquare$

Consider symplectic disk $\Sigma_-$ through the point $q(-N).$ In
normalized coordinates it has a representation $x=0,$ $y=\Psi(u,v).$ Hence
the tangent plane to it at $q(-N)$ has a representation $\xi=0,$ $\eta = a\chi^1 + b
\chi^2.$ This implies this plane to intersect transversely any Lagrangian
cylinder $C_n^-(I)$ and these intersections form a foliation of the plane
into ellipses. Similar foliation into ellipses exists in the tangent plane
to $\Sigma_+$ through the point $q(N)$, it is generated by intersection of
this plane with cylinders $C_n^+(I)$. Differential of the global map $\mathcal G =
f^{2N}$ is a linear symplectic map and its restriction to $\Sigma_-$ is a
linear symplectic map $L_N:T_{q(-N)}\Sigma_- \to T_{q(N)}\Sigma_+$. If we fix $I$, then
related ellipses $E_-(I),$ $E_+(I)$ in $\Sigma_-, \Sigma_+$, respectively, have the same area
and $L_N(E_-(I))$ intersects $E_+(I)$ at four different points
transversely, due to Assumption \ref{scattering_as}.

Let us fix a sufficiently small a neighborhood of fixed point $p$ on the center manifold
$W^c$. Its smallness is controlled by a parameter $\eps.$
To this end we, following \cite{Moser}, we introduce symplectic polar coordinates on $W^c = \{x=y=0\}$ in the neighborhood of fixed point
$$
u = \eps \sqrt{2I}\cos\theta,\quad v = \eps \sqrt{2I}\sin\theta.
$$
Then we get a symplectic map on $W^c$:
\begin{equation}\label{moser_variables}
\begin{array}{ccl}
\bar{I} & = & I + \mathcal{O}(\eps^3),\\
\bar{\theta} & = & \theta + \alpha + \nu \eps^2 I + \mathcal{O}(\eps^3).
\end{array}
\end{equation}
According to the Moser theorem \cite{Moser}, there is an $I_0>0$ such that for $\eps$ small enough and any given $I,$ $0<I<I_0$, such that
the number $\alpha + 2\nu\eps^2 I$ is Diophantine,
there exists an invariant curve which is $\eps^2$-close in $C^2$ topology in the space of curves $I=r(\theta)$ to curve $u^2+v^2 =2\eps^2 I$.
The map (\ref{moser_variables})
is the restriction of the initial map on a neighborhood of $p$ in $W^c$ and normalized up to the third order terms. Thus each invariant curve
on the center manifold has two invariant Lagrangian cylinders being its stable and unstable manifolds, they are $C^2$-close to
cylinders $u^2+v^2 =2\eps^2 I$, $y=0$, or $x=0$, respectively, of the truncated map.

Let us verify that traces of $W^s(\gamma)$, $W^u(\gamma)$, corresponding to invariant KAM curve $\gamma$ on
$T_p W^c$ with action $2\eps^2 I$ also intersect on $\Sigma_+$ transversally
along four points. To this purpose, we consider the restriction of $\mathcal{G}$ on $\Sigma_-$
near point $q(-N)$ with values on disk $\Sigma_+$ near point $q(N)$. Fix in $V$ some neighborhood of $p$
on $W^c$ defined by $\eps$ small enough,
and let $\gamma$ be some KAM-curve in this neighborhood. This defines some
$I$. The restriction of $\mathcal{G}$ on $\Sigma_-$ is a two dimensional
symplectic map (\ref{restr}). Since the coordinates on $\Sigma_-$ and $\Sigma_+$ are
$(u,v)$, after the change of variables  (\ref{moser_variables}) where $\eps$ and $I$ are considered
as parameters, we come to the system for intersection points of $W^s(\gamma),$ $W^s(\gamma)$. Taking into account
that these manifolds of the same KAM curve, we get the value $I$ be the same and then we have:
$$
\begin{array}{ccl}
\eps\sqrt{2I}\cos\bar{\theta} & = & \eps a_{11}\sqrt{2I}\cos\theta +
                                    \eps a_{12}\sqrt{2I}\sin{\theta} + \mathcal{O}(\eps^2),\\[\jot]
\eps\sqrt{2I}\sin\bar{\theta} & = & \eps a_{21}\sqrt{2I}\cos\theta +
                                    \eps a_{22}\sqrt{2I}\sin{\theta} + \mathcal{O}(\eps^2).
\end{array}
$$
Dividing the equations on $\eps\sqrt{2I}$, squaring both sides of each of the equalities and sum them we get the equation for $\theta$
corresponding intersection points:
$$
1 = (a_{11}\cos\theta + a_{12}\sin{\theta})^2 + (a_{21}\cos\theta + a_{22}\sin{\theta})^2 + \mathcal{O}(\eps).
$$
This equation have precisely four simple roots, if matrix $A$ is not a rotation matrix \cite{Le2}.
Therefore, due to implicit function theorem,
traces of cylinders of truncated map on disk $\Sigma$ intersect transversally along four points.
This implies that $C^2$-close traces of cylinders of full map,
that is curves $w_u(\gamma)$ and $w_s(\gamma)$, also intersect transversally along four points. Theorem \ref{main_th}
has been proved.$\blacksquare$


\section{Homoclinic orbits to saddle periodic orbits on $W^c$}
As is well known, a sufficiently smooth 2-dimensional symplectic diffeomorphism near its generic elliptic fixed
point $O$ (some Birkhoff coefficient in the normal form does not vanish) is a twist map with respect to
action-angle variables given by symplectic coordinates of the Birkhoff normal form. In particular, this implies the
existence of a positive measure Cantor set of invariant KAM curves with Diophantine rotation numbers
accumulating at $p$ \cite{Moser}. For our case as such diffeomorphism we have the restriction of $f$ onto
$W^c(p)$. Since $W^c$ is a smooth normally hyperbolic invariant
submanifold for $f$ the Fenichel results \cite{Fenichel1,Fenichel2} apply
that gives a smooth foliation of $W^{cs}$ and $W^{cu}$ into smooth curves (see Appendix).
This smooth foliations allow us to define two smooth maps $F_u: W^c \to
\Sigma$, $F_s: W^c \to \Sigma$.

Let us fix some such invariant KAM curve $\gamma$ with a Diophantine rotation number. Near such the curve there is a positive
measure set of other smooth invariant curves accumulating at $\gamma$ in at least $C^2$-topology. Another consequence of the KAM theorem,
twist condition is that for any rational number $p/q$ with incommensurate integers $p,q$ there are at least two $q$-periodic points.
Generically, these two periodic orbits are elliptic one and hyperbolic another. Moreover, if this diffeomorphism
satisfies some additional genericity condition (sometimes, it is called as the Birkhoff genericity \cite{Robin}),
then stable and unstable separatrices of the hyperbolic orbit
intersect transversely along related homoclinic orbits in $W^c$. This
allows, in particular, to construct a ``fence'' made up of stable and unstable
separatrices separated one invariant curve from another one.

\begin{figure}
\centering
\parindent=0pt
\includegraphics[width=1.2\textwidth]{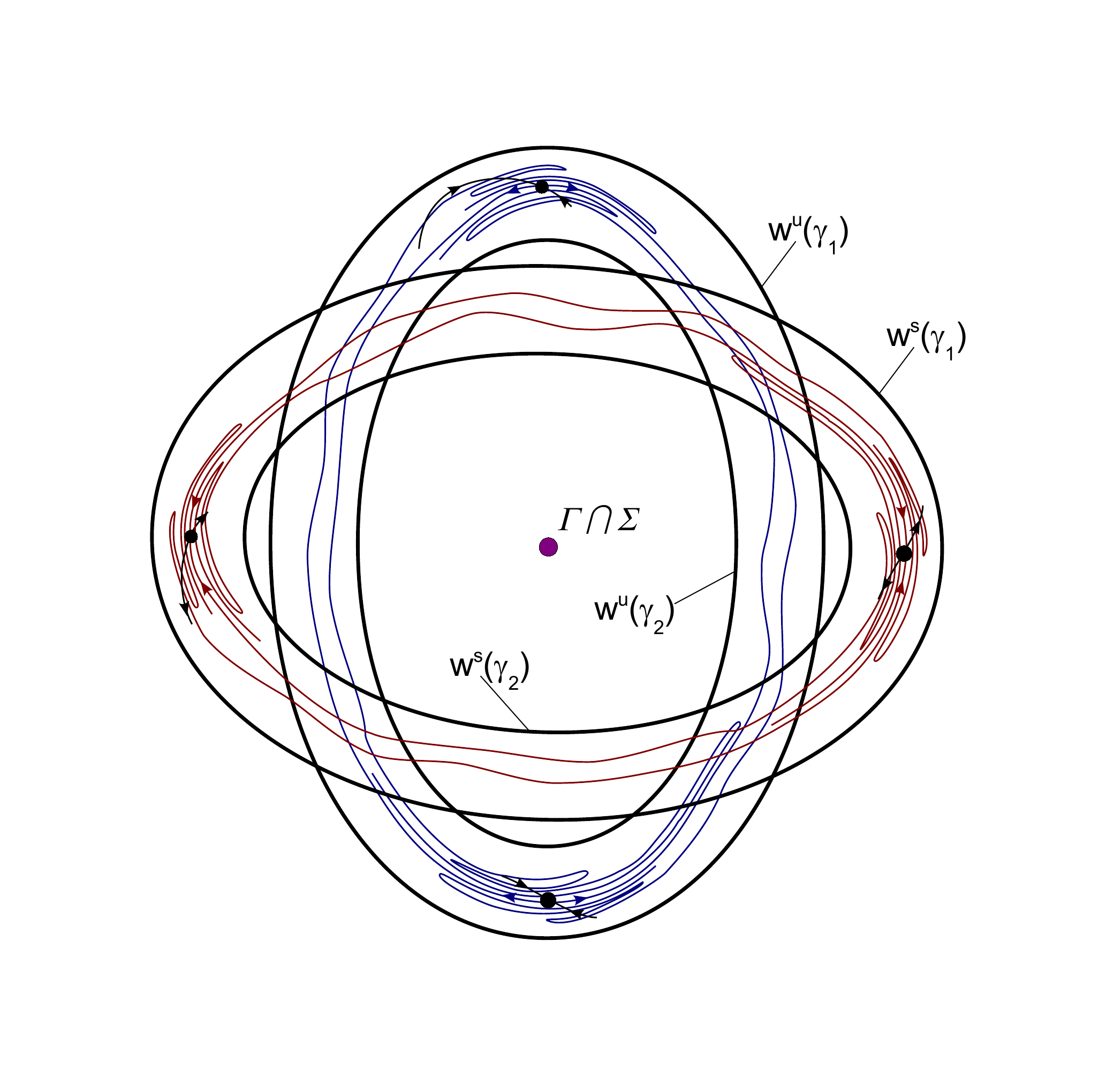}
\caption{A fence on $\Sigma$ made up of stable and unstable manifolds of a saddle periodic orbit.}
\label{fig:fence}
\end{figure}

This can be done in the following way. Let us assume $q=2$ to be definite. Take
one hyperbolic 2-periodic point $m_0$ and let $m_1$ be its first iteration: $m_1
=f(m_0).$ Suppose that unstable manifold $W^u(m_0)$ transversely intersects
stable manifold of $W^s(m_1)$ and stable manifold of $W^s(m_0)$ transversely intersects
unstable manifold of $W^u(m_1)$. In this case, due to so-called lambda lemma
\cite{Smale}, the topological limit of $W^u(m_0)$ contains
$W^u(m_1)$ and vise versa. Thus we get some closed invariant set $\mathfrak F_u$ made up of these curves and their closures.
Similar set $\mathfrak F_s$ is formed by stable manifolds (see Fig.\ref{fence}).

Take $\mathfrak F_u$ and for every its point consider the related unstable
leaf of the unstable foliation in $W^{cu}$. Then map $F_u$ transforms set $\mathfrak
F_u$ to the homeomorphic set in $\Sigma.$ Similar set in $\Sigma$ is obtained from $\mathfrak
F_s$ using $F_s$. Now choose any invariant KAM curve in $W^c$ in a neighborhood
$U$ where Theorem \ref{main_th} applies. Choose a sufficiently close
invariant KAM curve $\gamma_1$ such that on $\Sigma$ related traces $w^s(\gamma),w^s(\gamma_1)$
intersect transversely traces $w^u(\gamma),w^u(\gamma_1)$. Then we have on
$\Sigma$ two annuli: $A_s$ bounded by $w^s(\gamma),w^s(\gamma_1)$ and
$A_u$ bounded by $w^u(\gamma),w^u(\gamma_1)$. These annuli intersect each
other in such a way that each boundary curve of one annulus intersects
every boundary curve of another annulus transversely. Since
the restriction $f_c$ of $f$ on $W^c$ is a twist map, then invariant KAM curves
$\gamma, \gamma_1$ have different rotation numbers $\rho, \rho_1$. Thus
there are periodic orbits inside the annulus between $\gamma, \gamma_1$
corresponding to some rational $\rho < \rho_* < \rho_1,$ $\rho_* = p/q.$ If $f_c$ is
Birkhoff generic, then the half of these periodic orbits are hyperbolic
Birkhoff $q$-periodic and its stable manifolds form a fence $\mathfrak F_s$ in
$A_s \subset \Sigma$. This fence separates $A_s$ in the sense that if we
take two points on different boundary curves of $A_s$, then any path going
from one point to another one will cut $\mathfrak F_s$. The same holds for
$\mathfrak F_u$ in $A_u.$ This implies
\begin{Th}\label{hom_Poin}
The sets $\mathfrak F_s$, $\mathfrak F_u$ intersect, hence there are
Poincar\'e homoclinic orbits to a saddle hyperbolic periodic orbit on $W^s.$
\end{Th}
It is clear that in fact there are countably many such Poincar\'e homoclinic
orbits. It is impossible to assert that they are transverse or tangent
since this cannot be caught by such considerations.
\section{Multidimensional extension}

The problem we have studied possesses a multidimensional extension.
The tool to get this extension are essentially the same, so
we present only the related set up and formulations.
In a smooth symplectic manifold $(M, \Omega)$ of dimension $2n+2$ we
consider a symplectic diffeomorphism $f$ that possesses a fixed point $p$
of the elliptic-hyperbolic $(2,2n)$-type. The latter means this point has
the linearization operator with the only pair of complex eigenvalues on the unit
circle $\exp[\pm i \alpha]$ and remaining $2n$ eigenvalues are off the unit
circle and thus are met either in real pairs $\mu_k, \mu_k^{-1}$, $|\mu_k|<
1,$ or in complex quartets $\rho_m\exp[\pm i\kappa_m],$ $\rho^{-1}_m\exp[\pm
i\kappa_m],$, $\kappa_m \ne 0,\pi.$ Here one has $k + 2m = n.$ Such a
fixed point has locally a smooth two dimensional center manifold $W^c$
corresponding to the pair $\exp[\pm i \alpha]$ on which $p$ is an elliptic
fixed point and we assume henceforth it to be of generic elliptic type. Besides center
manifold, through fixed point other smooth manifolds pass: $n$-dimensional strong
stable $W^s$ and strong unstable $W^u$ ones, as well as $(n+2)$-dimensional center stable
$W^{cs}$ and center unstable $W^{cu}$ ones.

Analogs of three Assumptions 1-3 are

\begin{As}[Homoclinic intersection]\label{m-intersection_as}
Manifolds $W^u(p)$ and $W^s(p)$ have an intersection at some point $q$, generating thus a homoclinic orbit $\Gamma$ to point $p$.
\end{As}
Dimensions of stable $W^s$ ($n$) and center unstable manifold $W^{cu}$ ($n+2$) are complementary, this
allows one to assume
\begin{As}[Transversality condition]\label{m-transversality_as}
The intersection of manifolds $W^s(p)$ and $W^{cu}(p)$ at point $q$ is transverse.
\end{As}
Below we shall show that the linearized along the homoclinic orbit $\Gamma$ the
sequence of linearized map generates the linear symplectic scattering map $S: T_pW^c \to
T_pW^c$. We assume this map being generic that means as above that the
foliation into ellipses on the tangent plane $T_pW^c$ generated by the
linearized map $Df_p$ has the property: any ellipse $E$ of this foliation
satisfies
\begin{As}[Genericity condition]\label{m-scattering_as}
The intersection $E \cap S(E)$ is transverse and consists of four points.
\end{As}
If these three conditions hold then the analog of the main theorem is valid.
\begin{Th}\label{anmain_th}
Let a $C^r$-smooth, $r\ge 6$, symplectic diffeomorphism $f$ on a $C^\infty$-smooth $2(n+1)$-dimensional symplectic manifold $M$ with an
elliptic-hyperbolic fixed point of the type $(2,2n)$ obeys Assumptions \ref{m-intersection_as},
\ref{m-transversality_as}, \ref{m-scattering_as}. Then there is a neighborhood $U$ of homoclinic orbit $\Gamma$
such that every closed invariant KAM-curve on $W^c$ possesses four transverse homoclinic orbits in $U$.
\end{Th}
To prove this theorem we again first study the linearized nonautonomous
problem given by the linearization of $f$ on the homoclinic orbit
$\Gamma.$ Also, in order to avoid possible complications, one assumes in
addition that orbit $\Gamma$ leaves from $p$ and enters to $p$ along
leading direction in $W^u$ and $W^s$ (one or two dimensional).

Then, as above, we construct scattering map $S$ acting on $T_pW^c$ and assuming Assumptions
4-6 to hold we prove the Theorem. This proof uses again that the
transverse intersection of $W^{cu}$ and $W^{cs}$ near a homoclinic point $q \in
\Gamma$ occurs along a 2-dimensional disk $\Sigma$ which belongs to both
of them. Hence $(n+1)$-dimensional stable and unstable manifolds of any invariant
KAM curve $\gamma \subset W^c$ when continuing by $f$ in $W^{cs}$, $W^{cs}$, respectively, intersect again $\Sigma$ along closed
curves $w^s(\gamma)$, $w^u(\gamma)$. Genericity Assumption
\ref{m-scattering_as} implies this intersection to happen transversely at
four points through which homoclinic orbits to $\gamma$ pass.
\section{Appendices}

\subsection{Straightening invariant manifolds} \label{straightened_coordinates}

In some neighborhood of the fixed point the symplectic diffeomorphism under consideration can be written as
(\ref{linearized_map}). In this form 1-dimensional stable manifold $W^s$ is given as a smooth curve tangent to
the $x$-axis (at point $p$) and 1-dimensional unstable manifold $W^u$ is given as a smooth curve tangent to the
$y$-axis. Center stable $W^{cs}$ and center unstable $W^{cu}$ manifolds are given as graphs of the functions
$y=F(x,u,v)$ and  $x=G(y,u,v)$ being tangent at $p$ to 3-dimensional planes $y=0$ and $x=0$, respectively. Let us
first straighten the curves $W^s$, $W^u$:

\begin{Lm}
Let in linear symplectic space $(\mathbb R^4, \Omega=dx \wedge dy + du \wedge dv)$ a smooth curve $(x,y(x),u(x),v(x))$
through the point $(0,0,0,0)$ is given, such that $y'(0)=u'(0)=v'(0)=0$. Then this curve can be transformed by a
symplectic transformation $(x,y,u,v)\to (\xi,\eta,\nu,\omega)$ into the $\xi$-axis.
\end{Lm}
There are many of such transformations, for instance, this is one of them:
$$
\left\{
\begin{array}{ccl}
\xi & = & x \\
\eta & = & y-y(x)-v'(x)(u-u(x))+u'(x)(v-v(x)) \\
\nu & = & u-u(x) \\
\omega & = & v-v(x).
\end{array}
\right.
$$
All other transformations in a neighborhood of $p$ we perform holding straight $W^s$ and $W^u$. At the next step
we straighten $W^{cs}$ and $W^{cu}$:

\begin{Lm}
In some neighborhood of point $p$ there exist symplectic coordinates $(x,y,u,v)$, $\Omega = dx \wedge dy + du \wedge dv$,
such that submanifolds $W^{cs}, W^{cu}$ in these coordinates become flat, that is they are given as $y=0$, $x=0$, respectively.
\end{Lm}
{\bf Proof}. In principal, this lemma follows from the related result of the theory of symplectic manifolds (the relative
Darboux theorem) \cite{AG}. But for the reader's convenience, we present a direct proof. We follow the lines of the proof
of the Darboux theorem given in \cite{Arn}.

In coordinates (\ref{linearized_map}) $W^{cs}$ is expressed as $y=F(x,u,v)$, where
$F(0,0,0) = F_x(0,0,0) = F_u(0,0,0) = F_v(0,0,0) = 0$. Center unstable manifold $W^{cu}$ in the same coordinates is given
as $x = G(y,u,v)$, $G(0,0,0) = G_y(0,0,0) = G_u(0,0,0) = G_v(0,0,0) = 0$. Take function $p_1=y-F$ as a Hamilton function
and consider the related Hamiltonian flow. Since $\dot x = 1$ then in a neighborhood of $p$ small enough submanifold
$x=G(y,u,v)$ is transversal to flow orbits. We take this manifold as a cross-section to the flow. Denote $q_1(x,y,u,v)$
the time needed for the flow orbit through the initial point $(x_0 = G(y_0,u_0,v_0),y_0,u_0,v_0)$ to reach the point
$(x,y,u,v)$. Then $q_1=0$ for points on the cross-section and $p_1=0$ on $W^{cs}$ since it is a level of the Hamiltonian.
The Lie derivative of $q_1$ w.\,r.\,t. the vector field $X_{p_1}$ is equal to 1. Therefore, Hamiltonian vector fields with
the Hamilton functions $p_1, q_1$ are commute and independent in a neighborhood of $p$. Thus, orbits of $\mathbb R^2$-action
generated by these two commuting Hamilton functions give a smooth foliation into 2-dimensional orbits near $p$ and its
leaves are transversal to 2-dimensional submanifold being joint level of functions $p_1$ and $q_1$. Next we take joint
level $p_1=q_1=0$ that is just locally $W^c$. We introduce any local symplectic w.\,r.\,t. the restriction of 2-form
$\Omega$ on $W^c$ coordinates $(p_2,q_2)$ near $p$. These coordinates are extended onto a neighborhood of $W^c$ setting
 $(p_2, q_2)$ constant along the whole 2-dimensional orbit of the action through the point on $W^c$ with coordinates
$(p_2, q_2)$ on it. $\blacksquare$

\begin{Rm}
If $W^s$ and $W^u$ were previously made straighten, then one has $F(x,0,0) \equiv 0$, $G(y,0,0)\equiv 0$ and straightening
$W^{cs}, W^{cu}$ preserves $W^s, W^{u}$ straighten.
\end{Rm}

\subsection{Normal form near 1-elliptic fixed point}

Here we shall derive the normal form up to the terms of third order for a smooth symplectic 4-dim diffeomorphism
$f$ in neighborhood $V$ of its fixed 1-elliptic point $p$. Without a loss of generality one may assume $\alpha \in (0,\pi)$.
\begin{Prop}
In some neighborhood $V$ of fixed 1-elliptic point $p$ there exist symplectic coordinates $(x,y,u,v)$,
$\Omega = dx\wedge dy + du \wedge dv$, such that diffeomorphism $f$ has the following form in these coordinates:
\begin{equation}\label{normal_diffeomorphism}
\left\{
\begin{array}{ccl}
x_{n+1} & = & \mu x_n(1 + a x_n y_n + b(u_n^2 + v_n^2)+ \mathcal{O}_3)\\[\jot]
y_{n+1} & = & \mu^{-1} y_n(1 - a xy - b(u_n^2+v_n^2)+ \mathcal{O}_3) \\[\jot]
u_{n+1} & = & u_n \cos(\alpha + \nu(u_n^2 + v_n^2)) - v_n \sin(\alpha + \nu(u_n^2 + v_n^2)) \\[\jot]
        &   & \multicolumn{1}{r}{-\kappa x_n y_n (u_n\sin\alpha + v_n\cos\alpha) + \varphi(x_n,y_n,u_n,v_n)}\\[\jot]
v_{n+1} & = & u_n \sin(\alpha + \nu(u_n^2 + v_n^2)) + v_n \cos(\alpha + \nu(u_n^2 + v_n^2)) \\[\jot]
        &   & \multicolumn{1}{r}{ + \kappa x_n y_n (u_n \cos\alpha - v_n \sin\alpha) + \psi(x_n,y_n,u_n,v_n),}
\end{array}
\right.
\end{equation}
where $\nu \ne 0$, functions $\varphi,\,\psi$ are of the fourth order and higher at the origin,
$\mathcal{O}_3$ means terms of third order and higher at the origin. In these coordinates manifolds
$W^s$, $W^u$ coincide with $x$-axis, $y$-axis, respectively, that is the following identities hold:
$\varphi(x,0,0,0)=\varphi(0,y,0,0)=\psi(x,0,0,0)=\psi(0,y,0,0)\equiv 0$.
\end{Prop}
{\bf Proof}. At the first step we straighten manifolds $W^{cs}, W^{cu}$ in the neighborhood of
$p$ (see Appendix). As the result, first two relations in (\ref{linearized_map}) are transformed to the form:
$$
x_{n+1} = x_n(\mu + \cdots),\quad
y_{n+1} = y_n(\mu^{-1} + \cdots).
$$
Next we apply the standard normal form method for symplectic maps (see for instance \cite{AG}). We shall use such symplectic coordinate
transformations which hold $W^{cs}$, $W^{cu}$ be straightened. Next we use complex coordinates instead of
$u,v$ in order to diagonalize the linear part of the third and fourth relations. Monomials of the second
order and those of third order other than resonance monomials can be killed.

Resonance relations for the set of eigenvalues $(\mu,\mu^{-1},e^{i\alpha},e^{-i\alpha})$ and integer
vectors $(m_1,m_2,n_1,n_2)$ have the form:
$$
\begin{array}{ccc}
\mu^{m_1-m_2-1}e^{i\alpha(n_1-n_2)} & = & 1,\\
\mu^{m_1-m_2+1}e^{i\alpha(n_1-n_2)} & = & 1,\\
\mu^{m_1-m_2}e^{i\alpha(n_1-n_2-1)} & = & 1,\\
\mu^{m_1-m_2}e^{i\alpha(n_1-n_2+1)} & =& 1.
\end{array}
$$
These relations can be rewritten in the following way:
\begin{equation}\label{resonance_1}
m_1 = m_2+1,\; \alpha(n_1-n_2) = 2\pi k, k\in \Z;
\end{equation}
\begin{equation}\label{resonance_2}
m_2 = m_1+1,\; \alpha(n_1-n_2) = 2\pi k, k\in \Z;
\end{equation}
\begin{equation}\label{resonance_3}
m_1=m_2,\; \alpha(n_1-n_2-1) = 2\pi k, k\in \Z;
\end{equation}
\begin{equation}\label{resonance_4}
m_1=m_2,\; \alpha(n_1-n_2+1) = 2\pi k, k\in \Z.
\end{equation}
From these relations for integers vectors $(m_1,m_2,n_1,n_2)$ such that $|m|+|n|= m_1 + m_2 + n_1 + n_2 = 2 ,\, 3$
we derive that if $|m|+|n|=2,$ then resonance relations (\ref{resonance_1}) and (\ref{resonance_2}) are absent,
resonance relations (\ref{resonance_3}) and (\ref{resonance_4}) are the same as for the case of 2-dimensional
elliptic point: $\alpha = 2\pi/3$. Thus, according to our assumptions (that $p$ is 1-elliptic fixed point of
generic type for $f|_{W^c}$) resonant monomials of second order can be removed. If $|m|+|n|=3$, then for
resonance relations (\ref{resonance_1}) we get resonant monomials $x^2y$, corresponding to 2-dimensional
saddle point of symplectic diffeomorphism, and $x(u^2+v^2)$. For relations (\ref{resonance_2}) we get
$xy^2$ and $y(u^2+v^2)$. Under an assumption that strong resonances are absent in the system, relations
(\ref{resonance_3}) and (\ref{resonance_4}) show resonant monomials $xyu$, $xyv$, $u(u^2+v^2)$ and
$v(u^2+v^2)$ cannot be eliminated. The last two present in the normal form for a diffeomorphism in
a neighborhood of an elliptic point. Taking into account that the transformation to the normal form
should be symplectic we get (\ref{normal_diffeomorphism}).
$\blacksquare$

\subsection{Invariant foliations and their smoothness}

In this subsection we verify the conditions from \cite{Fenichel1,Fenichel2} which guarantee the existence and smoothness of
stable and unstable invariant foliations within manifolds $W^{cs},$ $W^{cu}$, respectively. In particular, these conditions
imply the existence of stable and unstable smooth invariant manifolds for KAM-curves on the center manifold $W^c$.
Homoclinic orbits to KAM-curves belongs to intersection of these manifolds. Note that this fact does not follow
immediately from the Hirsch-Pugh-Shub theorem \cite{HPS} and we use the theory developed by Fenichel \cite{Fenichel1,Fenichel2}.
Let us recall the definition of weakly overflowing invariant set according to \cite{Fenichel1}:

\begin{Def}[N.\,Fenichel, 1974]
Let $U$ and $V$ be open subsets of some $C^l$-manifold $M_1$, $1 \leq l \leq \infty$, and let $F:U\rightarrow V$ be a
$C^l$-diffeomorphism. A set $\Lambda \subset U$ is called weakly overflowing invariant $($under $F$$)$ if
$\Lambda \subset F(\Lambda)$.\\[2\jot]
Let $TF: TU \rightarrow TV$ be the map induced by $F$ on tangent spaces. A sub-bundle $E \subset TM_1 \big|_{\Lambda}$
is called weakly overflowing invariant if $E \subset TF(E)$.
\end{Def}

Let us choose any invariant KAM-curve on $W^c$ in a sufficiently small neighborhood of $p$. Then the closure of subset
in $W^c$ bounded by this KAM-curve is weakly overflowing invariant set under diffeomorphism $f$. Denote this compact
set as $\Lambda$. Here we assume $M_1=W^{cu}$ and consider $f\big|_{M_1}$. To apply expanding family theorems one
need to show that there exists weakly overflowing invariant sub-bundle $\textsl{E} \subset TM_1 \big|_\Lambda$.
It will be proved using contraction mapping principle. Recall that locally near $p$ manifolds $W^{cu}, W^{cs}$
are straightened (i.e. $x\equiv 0$ on $W^{cu}$ and $y\equiv 0$ on $W^{cs}$) and $(y,u,v)$ are coordinates on it.
Then the restriction $f$ on $W^{cu}$ has following form:
$$
\left\{
\begin{array}{ccl}
    y_{n+1} & = & \mu^{-1}y_n + h(y_n,u_n,v_n)y_n\\
    u_{n+1} & = & u_n \cos \alpha - v_n \sin \alpha + g_1(y_n,u_n,v_n)\\
    v_{n+1} & = & u_n \sin \alpha + v_n \cos \alpha + g_2(y_n,u_n,v_n),
\end{array}
\right.
$$
where functions $g_1, g_2$ are of second order in $(y,u,v)$, $h$ is first order function. Now let us change
variables $(y_n,u_n,v_n) \rightarrow (\hat{y}_n,\hat{u}_n,\hat{v}_n)$:
$$
\left\{
\begin{array}{ccl}
y_n & = & \hat{y}_n \\
\left(
  \begin{array}{c}
    u_n \\ v_n
  \end{array}
\right) & = & R_{n\alpha}
\left(
  \begin{array}{c}
    \hat{u}_n \\ \hat{v}_n
  \end{array}
\right).
\end{array}
\right.
$$
In new coordinates the diffeomorphism will have the form (we keep old notations for variables):
\begin{equation}\label{diffeomorphism_kam_curve_manifold}
\left\{
\begin{array}{ccl}
    y_{n+1} & = & \mu^{-1}y_n + h(y_n,u_n,v_n)y_n\\
    u_{n+1} & = & u_n + g_1(y_n,u_n,v_n)\\
    v_{n+1} & = & v_n + g_2(y_n,u_n,v_n),
\end{array}
\right.
\end{equation}
At any point $m \in \Lambda$ differential $Df$ has the following representation (recall $\Lambda \subset W^c$ is given as $y \equiv 0$):
$$
Df_{m} =
\left(
    \begin{array}{ccc}
            \mu^{-1} + h(m) & 0 & 0 \\
            g_{1_{y}}(m) & 1 + g_{1_{u}}(m) & g_{1_{v}}(m) \\[\jot]
            g_{2_{y}}(m) & g_{2_{u}}(m) & 1 + g_{2_{v}}(m)
    \end{array}
\right).
$$
Since we work in one coordinate chart $(y,u,v)$, we will denote $(\eta, \chi_1, \chi_2)$ as coordinates
in the tangent space to a point $m$.
Now consider any orbit $(\ldots, m_n, m_{n+1}, \ldots)$ of $f\big|_{W^c}$, $m=m_0,$ which belongs to
$\Lambda$ (that is to $W^c$). At each point of this orbit choose in $T_{m_n}M_1$ a straight line through
the origin in the tangent space being transversal to plane $\eta=0$. Such straight line can be given
parametrically: $(\eta, p_n\eta, q_n\eta)$, $\eta \in \R$, functions $p_n$ and $q_n$ smoothly depend
on $m$ in $\Lambda$. Differential $Df_{m}$ transforms this line to another line in $T_{f(m)}M_1$:
\begin{equation}\label{differential_map}
\begin{array}{ccl}
\bar{\eta} & = & (\mu^{-1} + h(m_n))\eta, \\[2\jot]
\bar{p}_{n+1} & = & \frac{\displaystyle g_{1_{y}}(m_n) + (1+g_{1_{u}}(m_n))p_n + g_{1_{v}}(m_n)q_n}
                         {\displaystyle \mu^{-1} + h(m_n)},\\[3\jot]
\bar{q}_{n+1} & = & \frac{\displaystyle g_{2_{y}}(m_n) + g_{2_{u}}(m_n)p_n + (1+g_{2_{v }}(m_n))q_n}
                         {\displaystyle \mu^{-1} + h(m_n)}.
\end{array}
\end{equation}
R.h.s. of the second and third relations define an operator in the Banach space of uniformly bounded sequences
$\zeta = \left\{(p_n,q_n)\right\}$ with norm $\left\| \zeta \right\| = \underset{n \in \Z}{\sup}
(\left|p_n\right|, \left|q_n\right|)$. Indeed, operator will transforms a uniformly bounded sequence
$\zeta$ to the bounded one, as all derivatives $(g_i)_{y}'$, $(g_i)_{u}'$, $(g_i)_{v}'$, $i = 1, 2$,
calculated at $m_n$, are small (functions $g_i$, $i = 1, 2$, are of second order at zero):
$$
\begin{array}{ccl}
\left| \bar{p}_{n+1} \right|& \leq &
\frac{\displaystyle \left|g_{1_{y}}(m_n)\right| + \left|1+g_{1_{u}}(m_n)\right|\left\|\zeta\right\|
                    + \left|g_{1_{v}}(m_n)\right|\left\|\zeta\right\|}
     {\displaystyle \left|\mu^{-1} + h(m_n)\right|},\\[4\jot]
\left| \bar{q}_{n+1} \right| & \leq &
\frac{\displaystyle \left|g_{2_{y}}(m_n)\right| + \left|g_{2_{u}}(m_n)\right|\left\|\zeta\right\|
                    + \left|1+g_{2_{v}}(m_n)\right|\left\|\zeta\right\|}
     {\displaystyle \left|\mu^{-1} + h(m_n)\right|}.
\end{array}
$$
To prove the operator is contracting, consider following inequalities:
\begin{multline*}
\left| \bar{p}_{n+1}^1 - \bar{p}_{n+1}^2 \right|
        \leq \frac{\displaystyle \left|1+g_{1_{u}}(m_n)\right|\left|p_n^1 - p_n^2\right|
        + \left|g_{1_{v}}(m_n)\right|\left|q_n^1 - q_n^2\right|}
        {\displaystyle \left|\mu^{-1} + h(m_n)\right|}\\
\shoveright{\leq \frac{\displaystyle \left( \left|1+g_{1_{u}}(m_n)\right|+ \left|g_{1_{v}}(m_n)\right| \right)
        \left\|\zeta^1 - \zeta^2\right\|} {\displaystyle \left|\mu^{-1} + h(m_n)\right|},}\\
\shoveleft{\left| \bar{q}_{n+1}^1 - \bar{q}_{n+1}^2 \right|
        \leq \frac{\displaystyle \left|g_{2_{u}}(m_n)\right|\left|p_n^1 - p_n^2\right|
        + \left|1+g_{2_{v}}(m_n)\right|\left|q_n^1 - q_n^2\right|}
        {\displaystyle \left|\mu^{-1} + h(m_n)\right|}}\\
\leq \frac{\displaystyle \left( \left|g_{2_{u}}(m_n)\right| + \left|1+g_{2_{v}}(m_n)\right| \right)
         \left\|\zeta^1 - \zeta^1\right\|}
         {\displaystyle \left|\mu^{-1} + h(m_n)\right|}.
\end{multline*}
Quantities such as
$$
\frac{\displaystyle \left|1+g_{1_{u}}(m_n)\right|+ \left|g_{1_{v}}(m_n)\right|}
     {\displaystyle \left|\mu^{-1} + h(m_n)\right|}
$$
are less than 1 uniformly in $n$, if we are working in a neighborhood of $p$ small enough, so that
the operator is contracting. According to contraction mapping principle there exists a unique fixed
point of the operator, namely some sequence $\left\{(p^*_n,q^*_n)\right\}$. This sequence (as the
point in the related Banach space) depends continuously in $m\in \Lambda$. The straight lines corresponding
to this sequence, as $m$ varies along $\Lambda$, form weakly overflowing invariant sub-bundle
$\textsl{E} \subset TM_1 \big|_\Lambda$.

Next we choose a vector bundle $N \subset TM_1 \big|_\Lambda$, complementary to $E$, that is
$TM_1 \big|_\Lambda = E \oplus N$. We set $N = T\Lambda$. According to \cite{Fenichel1} for
$m \in \Lambda$, any $v_0 \in E$ and $W_0 \in N$ let
$$
\begin{array}{ccl}
v_{-k} & = & Df^{-k}(m)v_0, \\[\jot]
w_{-k} & = & \pi^{N} Df^{-k}(m)w_0,
\end{array}
$$
where $\pi ^{N}$ is projection to $N$ (note that in our case $N$ is invariant under $Df$, so one can
just let $w_{-k} = Df^{-k}(m)w_0$). Let us also define two numbers
$$
\begin{array}{ccl}
\alpha^*(m) & = & \underset{\alpha > 0}{\inf} \left\{\left| v_{-k} \right| / \alpha^k \rightarrow 0
                  \text{ as } k\rightarrow \infty, \forall\; v_0 \in E \right\}, \\[3\jot]
\rho^*(m) & = & \underset{\rho > 0}{\inf}
                \left\{\left| v_{-k} \right| / \left| w_{-k} \right| / \rho^k
                \rightarrow 0 \text{ as } k \rightarrow \infty,
                \forall\; v_0 \in E, w_0 \in N \right\}.
\end{array}
$$
The number $\alpha ^*(m)$ is an asymptotic measure of the growth of vectors in $E$ under the action of $Df^{-1}$,
and $\rho^*(m)$ is an asymptotic measure of the ratio of the growth of vectors in $E$ to the growth of vectors
in $N$ under the action of  $Df^{-1}$.

Recall one more definition and formulate expanding family theorem \cite{Fenichel1} for reader's convenience.
\begin{Def}[N.\,Fenichel, 1974]
The pair $(\Lambda, E)$ is called an invariant set with expanding structure for $f$, if $\Lambda$ is compact
and weakly overflowing invariant, $E$ is weakly overflowing invariant, and $\alpha^*(m)<1$, $\rho^*(m)<1$
for all $m \in \Lambda$.
\end{Def}
\begin{Th}[Expanding Family Theorem, N.\,Fenichel, 1974]
Let $M_1$ be $C^l$-manifold, $1 \leq l \leq \infty$, and let $F: M_1 \rightarrow M_1$ be a $C^l$-diffeomorphism.
Let $(\Lambda, E)$ be an invariant set with expanding structure. Then there is a family of $C^l$-manifolds
$W^E(m)$, $m \in \Lambda$, invariant in the sense that
$$
F^{-1}(W^E(m)) = W^E(F^{-1}(m)).
$$
The manifold $W^E(m)$ is $C^l$-diffeomorphic to the fiber $E_m$, and is tangent to $E_m$ at $m$.
\end{Th}

Let us show that in our case $\alpha^*(m) < 1$, $\rho^*(m) < 1$ for any $m \subset \Lambda$. Take any $v_0 \in E$,
it has the coordinate representation: $v_0 = (\eta_0, p^*_0 \eta_0, q^*_0 \eta_0)$. Vector $v_{-k}$ will have
the representation: $v_{-k} = (\eta_{-k}, p^*_{-k} \eta_{-k}, q^*_{-k} \eta_{-k})$. Taking into account first
equality from (\ref{differential_map}) one gets:
$$
\begin{array}{ccl}
\eta_{-1} & = & \frac{\displaystyle \mu \eta_0}{\displaystyle 1 + \mu h(m_{-1})},\\
\eta_{-2} & = & \frac{\displaystyle \mu \eta_{-1}}{\displaystyle 1 + \mu h(m_{-2})}
                = \frac{\displaystyle \mu^2 \eta_0}{\displaystyle (1 + \mu h(m_{-1}))(1 + \mu h(m_{-2}))},\\
\cdots    &   & \\
\eta_{-k} & = & \frac{\displaystyle \mu^k \eta_0}{\displaystyle (1 + \mu h(m_{-1})) \ldots (1 + \mu h(m_{-k}))}.
\end{array}
$$
Now consider the ratio
$$
\frac{\displaystyle \left| v_{-k} \right|}{\displaystyle \alpha^k}
= \frac{\displaystyle \left| \eta_{-k} \right| \sqrt{1 + (p^*_{-k})^2 + (q^*_{-k})^2 }}{\displaystyle \alpha^k}
= \frac{\displaystyle \mu^k \left| \eta_0 \right| \sqrt{1 + (p^*_{-k})^2 + (q^*_{-k})^2 }}
       {\displaystyle \alpha^k \left| (1 + \mu h(m_{-1})) \ldots (1 + \mu h(m_{-k})) \right|}.
$$
Quantity $\left| \eta_0 \right| \sqrt{1 + (p^*_{-k})^2 + (q^*_{-k})^2 }$ is bounded. Function $h$ is of the
first order, $m_n$ lies in small neighborhood of fixed point $p$. Let us define $\delta_1 =
\underset{n \in \Z}{\sup} |h(m_n)|$. This value is of order of size (radius) of neighborhood  and
hence is small enough. The following estimates are valid:
\begin{equation}\label{inequalities_for_alpha*}
\frac{\displaystyle \mu^k \left| \eta_0 \right| \sqrt{1 + (p^*_{-k})^2 + (q^*_{-k})^2 }}
     {\displaystyle \alpha^k (1 + \mu \delta_1)^k}
\leq
\frac{\displaystyle \left| v_{-k} \right|}{\displaystyle \alpha^k}
\leq
\frac{\displaystyle \mu^k \left| \eta_0 \right| \sqrt{1 + (p^*_{-k})^2 + (q^*_{-k})^2 }}
     {\displaystyle \alpha^k (1 - \mu \delta_1)^k}.
\end{equation}
Quantity in the r.\,h.\,s. of inequality (\ref{inequalities_for_alpha*}) tends to zero as
$k \rightarrow \infty$ if $\mu / \alpha (1 - \mu \delta_1) < 1$, that is $\alpha > \mu / (1 - \mu \delta_1)$.
On the other hand, quantity in the l.\,h.\,s. of inequality (\ref{inequalities_for_alpha*}) tends to zero as
$k \rightarrow \infty$ if $\alpha > \mu / (1 + \mu \delta_1)$. Thus, we get:
$$
\alpha^*(m_0) \in \left[ \frac{\displaystyle \mu}{\displaystyle 1 + \mu \delta_1},
                         \frac{\displaystyle \mu}{\displaystyle 1 - \mu \delta_1}
                  \right] < 1.
$$

The inverse map for $Df_{m_n}$ has following representation:
\begin{equation}\label{inverse_map_differential}
Df^{-1}_{m_n} =
\left(
    \begin{array}{ccc}
            \frac{\displaystyle \mu}{\displaystyle 1 + \mu h(m_n)} & 0 & 0 \\
            \ldots & 1 + l_1(m_n) & l_2(m_n) \\[\jot]
            \ldots & l_3(m_n) & 1 + l_4(m_n)
    \end{array}
\right),
\end{equation}
where dots, $l_1$, $l_2$, $l_3$ and $l_4$ are of at least first order functions. Take any $w_0 \in N$,
it has the form: $w_0 = (0,c_0,d_0)$, and let us consider $k$-th iteration of  $w_0$ under $DF^{-1}$:
$w_{-k} = (0,c_{-k},d_{-k})$. From (\ref{inverse_map_differential}) one gets that coordinates of
$w_{-k}$ change as follows:
$$
\begin{array}{ccl}
c_{-k} & = & (1 + l_1(m_{-k+1}))c_{-k+1} + l_2(m_{-k+1})d_{-k+1}\\
d_{-k} & = & l_3(m_{-k+1})c_{-k+1} + (1 + l_4(m_{-k+1}))d_{-k+1}.
\end{array}
$$
Denote $\delta_2 = \underset{i}{\sup}\, \underset{n \in \Z}{\sup} |l_i(m_n)|$, $i=\overline{1,4}$.
This quantity is also small enough and of the order of the size of neighborhood. Next estimates are valid:
$$
\begin{array}{ccl}
|c_{-k}| & \leq & (1 + \delta_2)|c_{-k+1}| + \delta_2|d_{-k+1}|
                  = |c_{-k+1}| + \delta_2(|c_{-k+1}| + |d_{-k+1}|)\\[\jot]
|d_{-k}| & \leq & \delta_2|c_{-k+1}| + (1 + \delta_2)|d_{-k+1}|
                  = |d_{-k+1}| + \delta_2(|c_{-k+1}| + |d_{-k+1}|),
\end{array}
$$
and, so:
$$
\begin{array}{ccl}
|c_{-k}| + |d_{-k}| & \leq & |c_{-k+1}| + |d_{-k+1}| + 2 \delta_2 (|c_{-k+1}| + |d_{-k+1}|)\\
                    & =    & (|c_{-k+1}| + |d_{-k+1}|)(1 + 2\delta_2)\\
                    & \leq & (|c_{-k+2}| + |d_{-k+2}|)(1 + 2\delta_2)^2 \leq \ldots\\
             \ldots & \leq & (|c_0| + |d_0|)(1 + 2\delta_2)^k.
\end{array}
$$
Then one gets:
\begin{multline*}
\sqrt{|c_{-k}|^2 + |d_{-k}|^2}
\leq \sqrt{2 \left( \max \left\{|c_{-k}|, |d_{-k}|\right\} \right)^2}
\leq \sqrt{2} \max \left\{|c_{-k}|, |d_{-k}|\right\}\\
\leq \sqrt{2} \left( |c_{-k}| + |d_{-k}| \right)
\leq \sqrt{2} (|c_0| + |d_0|)(1 + 2\delta_2)^k.
\end{multline*}
On the other side:
$$
\begin{array}{ccl}
|c_{-k}| & \geq & (1-\delta_2)|c_{-k+1}| - \delta_2|d_{-k+1}|
                  = |c_{-k+1}| - \delta_2(|c_{-k+1}| + |d_{-k+1}|)\\[\jot]
|d_{-k}| & \geq & (1 - \delta_2)|d_{-k+1}| - \delta_2|c_{-k+1}|
                  = |d_{-k+1}| - \delta_2(|c_{-k+1}| + |d_{-k+1}|),
\end{array}
$$
and, consequently,
$$
\begin{array}{ccl}
|c_{-k}| + |d_{-k}| & \geq & |c_{-k+1}| + |d_{-k+1}| - 2 \delta_2 (|c_{-k+1}| + |d_{-k+1}|)\\
                    & =    & (|c_{-k+1}| + |d_{-k+1}|)(1 - 2\delta_2)\\
                    & \geq & (|c_{-k+2}| + |d_{-k+2}|)(1 - 2\delta_2)^2 \geq \ldots\\
             \ldots & \geq & (|c_0| + |d_0|)(1 - 2\delta_2)^k.
\end{array}
$$
The following inequalities are valid:
$$
(|c_0| + |d_0|)(1 - 2\delta_2)^k \leq |c_{-k}| + |d_{-k}| \leq 2\max \left\{|c_{-k}|, |d_{-k}|\right\},
$$
so,
$$
\max \left\{|c_{-k}|, |d_{-k}|\right\} \geq \frac{\displaystyle |c_0| + |d_0|}{2} (1 - 2\delta_2)^k
$$
and
$$
\sqrt{|c_{-k}|^2 + |d_{-k}|^2}
\geq \sqrt{\left( \max \left\{|c_{-k}|, |d_{-k}|\right\} \right)^2}
= \max \left\{|c_{-k}|, |d_{-k}|\right\}\\
\geq \frac{\displaystyle |c_0| + |d_0|}{2}(1 + 2\delta_2)^k.
$$
Now let us evaluate
$$
\frac{\displaystyle \left| v_{-k} \right|}{\displaystyle \left| w_{-k} \right| \rho^k}
= \frac{\displaystyle \mu^k \left| \eta_0 \right| \sqrt{1 + (p^*_{-k})^2 + (q^*_{-k})^2 }}
       {\displaystyle \rho^k \left| (1 + \mu h(m_{-1})) \ldots (1 + \mu h(m_{-k})) \right|
                      \sqrt{|c_{-k}|^2 + |d_{-k}|^2}}.
$$
Next estimates are valid:
\begin{equation}\label{inequalities_for_rho*}
\begin{array}{ccl}
\frac{\displaystyle \left| v_{-k} \right|}{\displaystyle \left| w_{-k} \right| \rho^k}
& \geq &
\frac{\displaystyle \mu^k \left| \eta_0 \right| \sqrt{1 + (p^*_{-k})^2 + (q^*_{-k})^2 }}
     {\displaystyle \rho^k (1 + \mu \delta_1)^k \sqrt{2} (|c_0| + |d_0|) (1 + 2\delta_2)^k}\\[2\jot]
\frac{\displaystyle \left| v_{-k} \right|}{\displaystyle \left| w_{-k} \right| \rho^k}
& \leq &
\frac{\displaystyle 2 \mu^k \left| \eta_0 \right| \sqrt{1 + (p^*_{-k})^2 + (q^*_{-k})^2 }}
     {\displaystyle \rho^k (1 - \mu \delta_1)^k (|c_0| + |d_0|) (1 - 2\delta_2)^k}.
\end{array}
\end{equation}
The expression in the r.\,h.\,s. of the first inequality in (\ref{inequalities_for_rho*}) tends to
zero as $k$ tends to $\infty$ if $\rho > \mu /(1 + \mu \delta_1)(1 + 2 \delta_2)$, while expression
in the r.\,h.\,s. of the second inequality in (\ref{inequalities_for_rho*}) tends to zero as $k$ tends
to $\infty$ if $\rho > \mu / \displaystyle (1 - \mu \delta_1)(1 - 2 \delta_2)$. Thus, we get:
$$
\rho^*(m_0) \in \left[ \frac{\displaystyle \mu}{\displaystyle (1 + \mu \delta_1)(1 + 2 \delta_2)},\,
                       \frac{\displaystyle \mu}{\displaystyle (1 - \mu \delta_1)(1 - 2 \delta_2)}
                \right] < 1.
$$
Then expanding family theorem holds and for each point $m \in \Lambda$ there exists 1-dimensional
manifold (a curve) in $M_1$ being tangent to corresponding layer in $E$. The collection of these
curves is invariant under $f^{-1}$. In particular, collecting these manifolds for points of an
invariant KAM-curve defines its unstable manifold.

Now we want to have smoothness properties for the expanding foliation obtained. Let us apply smoothness
theorem for invariant sets with expanding structure to prove that these manifolds smoothly depend on
$m_0$ \cite{Fenichel2}. For this purpose we define the quantity $\tau^*$ according to \cite{Fenichel2}:
\begin{multline*}
\tau^*(m_0) = {\inf} \left\{ \tau:
              \left[ \left| v_{-k} \right| / \left| w_{-k} \right|
              \right]^{\tau} \left| \xi_{-k}\right| \rightarrow 0 \text{ as } k \rightarrow \infty, \right.\\
\left. \forall\; v_0 \in E, w_0 \in N, \xi_0 \in T_{m_0}\Lambda = N \right\}.
\end{multline*}

\begin{Th}[Smooth Invariant Bundle Theorem, N.\,Fenichel, 1977]
Let $U$ and $V$ be open subsets of a $C^l$-manifold $M_1$, and let $F: U \rightarrow V$ be a $C^l$-diffeomorphism,
$2 \leq l \leq \infty$. Let $\Lambda$ be a compact, properly embedded, $C^l$-manifold with boundary, overflowing
invariant under $F$. Let $(\Lambda,E)$ be an invariant set with expanding structure. If $1 \leq l' \leq l - 1$,
and $\tau^*(m) < 1/l'$ for all $m \in \Lambda$, then $E$ is a $C^{l'}$-smooth vector bundle.
\end{Th}

Taking into account estimates found before, one gets for diffeomorphism (\ref{diffeomorphism_kam_curve_manifold})
and any vector $\xi_0 = (0, \bar{c}_0, \bar{d}_0) \in N$, $\xi_{-k} = (0, \bar{c}_{-k}, \bar{d}_{-k})$, following
inequality to be valid:
\begin{multline}\label{inequality_tau*_1}
\frac{\displaystyle \left| v_{-k} \right|^{\tau}}
     {\displaystyle \left| w_{-k} \right|^{\tau}}|\xi_{-k}| \\
= \left( \frac{\displaystyle \mu^k \left| \eta_0 \right| \sqrt{1 + (p^*_{-k})^2 + (q^*_{-k})^2 }}
              {\displaystyle \left|(1 + \mu h(m_{-1})) \ldots (1 + \mu h(m_{-k})) \right|
                             \sqrt{|c_{-k}|^2 + |d_{-k}|^2}}
  \right)^{\tau} \sqrt{|\bar{c}_{-k}|^2 + |\bar{d}_{-k}|^2} \\
\leq \frac{\displaystyle \mu^{\tau k} C^{\tau}_1}
          {\displaystyle (1 - \mu \delta_1)^{\tau k} (1 - 2 \delta_2)^{\tau k}}
      \sqrt{2} (|\bar{c}_0| + |\bar{d}_0|) (1 + 2 \delta_2)^k,
\end{multline}
here $C_1$ is  constant which can be easily calculated. The r.\,h.\,s. of (\ref{inequality_tau*_1}) tends
to zero as $k \rightarrow \infty$ if
$$
\tau \geq \tau_1 =
\frac{\displaystyle \ln (1 + 2\delta_2)}{\displaystyle \ln [(1 - \mu \delta_1)(1 - 2\delta_2)] - \ln \mu} > 0.
$$
Quantity $\tau_1$ is small enough of order $\delta_2$. On the other hand,
\begin{equation}\label{inequality_tau*_2}
\frac{\displaystyle \left| v_{-k} \right|^{\tau}}
     {\displaystyle \left| w_{-k} \right|^{\tau}}|\xi_{-k}|
\geq \frac{\displaystyle \mu^{\tau k} C^{\tau}_2}
          {\displaystyle (1 + \mu \delta_1)^{\tau k} (1 + 2 \delta_2)^{\tau k}}
     \frac{|\bar{c}_0| + |\bar{d}_0|}{2} (1 - 2 \delta_2)^k.
\end{equation}
The r.\,h.\,s. of (\ref{inequality_tau*_2}) tends to zero as $k \rightarrow \infty$ if
$$
\tau \geq \tau_2 =
\frac{\displaystyle \ln (1 - 2\delta_2)}{\displaystyle \ln [(1 + \mu \delta_1)(1 + 2\delta_2)] - \ln \mu},\;
\tau_2 < 0.
$$
Therefore, $ \tau^*(m_0) \in \left[ \tau_2 , \tau_1 \right]$, so vector bundle $E$ is $C^{r'}$ smooth,
where $r' < 1 / \tau_1$, $r' \leq r - 2$.

Existence of stable manifold $W^s(\gamma)$ and its smoothness can be proved in a similar way.

\section{Acknowledgement}
The authors thank R. de la Llave and S.V. Gonchenko for useful discussions. We acknowledge a partial support from the
Russian Foundation for Basic Research under the grants 13-01-00589a (L.L.) and 14-01-00344
(A.M.). L.L. is also thankful for a support to the Russian Ministry of Science and Education
(project 1.1410.2014/K, target part) and the Russian Science Foundation (project
14-41-00044).


\end{document}